\documentclass{article}

\title{\LARGE \textbf{A Lower Bound for the Circumference Involving 
Connectivity\footnote{The original version is preprinted in Transactions
 of the Institute for Informatics and Automation Problems of the NAS 
(Republic of Armenia) and Yerevan State University, Mathematical Problems
 of Computer Science {\small 21 (2000) 129--155.}}
}}
\author{Zh.G. Nikoghosyan  \footnote{G.G. Nicoghossian up to 1997} \\
 \\Institute for Informatics and Automation Problems\\ National Academy of
 Sciences\\P. Sevak 1, Yerevan 0014, Armenia\\ E-mail: zhora@ipia.sci.am}

\begin{document}

\maketitle

\begin{abstract}
Let $G$ be a graph, $C$ a longest cycle in $G$ and $\overline{p}$, $\overline{c}$ the lengths of a longest 
path and a longest cycle in $G\backslash C$, respectively. Almost all lower bounds for the circumference base on
a standard procedure: choose an initial cycle $C_0$  in $G$ and try to enlarge it via structures of $G\backslash C_0$ and 
connections between $C_0$ and $G\backslash C_0$ closely related to 
$\overline{p}$, $\overline{c}$ and connectivity $\kappa$. Actually, each lower bound obtained in 
result of this procedure, somehow or is related to $\kappa$, $\overline{p}$, $\overline{c}$ but in forms of  
various particular values of $\kappa$, $\overline{p}$, $\overline{c}$ and the major
problem is to involve these invariants into such bounds as parameters. 
In this paper we present a lower
 bound for the circumference involving $\delta$, 
$\kappa$ and $\overline{c}$ and increasing with $\delta$, $\kappa$ and
 $\overline{c}$.
\end{abstract}

\section{Introduction}

We consider only finite undirected graphs without loops or multiple edges.
 A good reference for any undefined terms is [1]. The set of vertices of a
 graph $G$ is denoted by $V(G)$; the set of edges by $E(G)$. For $S$ a 
subset of $V(G)$, we denote by $G\backslash S$ the maximum subgraph of
 $G$ with vertex set $V(G)\backslash S$. For a subgraph $H$ of $G$ we 
use $G\backslash H$ short for $G\backslash V(H)$. 

Paths and cycles in a graph $G$ are considered as subgraphs of $G$. If 
$Q$ is a path or a cycle, then the length of $Q$, denoted by $|Q|$, 
is $|E(Q)|$. For $Q$ a path, we denote $|Q|=-1$ if and only if 
$V(Q)=\emptyset$. Throughout the paper, each vertex and edge can
be interpreted as cycles of lengths 1 and 2, respectively.

Let $G$ be a graph and $C$ a longest cycle in $G$. A cycle $C$ 
is a Hamilton cycle if $G\backslash C=\emptyset$ and is a dominating cycle
if $G\backslash C$ is edgeless. We use $n$ to denote the order,
$\delta$ the minimum degree and $\kappa$ the connectivity in $G$. The 
length $|C|$, denoted by $c$, is called a circumference. The lengths of a longest path and a 
longest cycle in $G\backslash C$, will be denoted by $\overline{p}$ and $\overline{c}$, respectively.

Almost all lower bounds for the circumference base on
a standard procedure: choose an initial cycle $C_0$  in a graph 
$G$ and try to enlarge it via structures of $G\backslash C_0$ and 
connections between $C_0$ and $G\backslash C_0$ closely related to 
connectivity $\kappa$ and $\overline{p}$, $\overline{c}$.     
Actually, each lower bound obtained in 
result of this procedure, somehow or is related to $\kappa$, $\overline{p}$, $\overline{c}$ but in forms of  
various particular values of $\kappa$, $\overline{p}$, $\overline{c}$ and the major
problem in long cycles theory is to involve these invariants into such bounds as parameters.

The starting result in this area, due to Dirac [2], bases on the minimum degree $\delta$.\\

\noindent\textbf{Theorem A} [2]. In every graph, $c\geq \delta+1$.\\

The second result in the same paper [2] shows that under 2-connectedness the
 bound $\delta+1$ in Theorem A can be replaced by $\min\{n,2\delta\}$. \\

\noindent\textbf{Theorem B} [2]. In every 2-connected graph, $c\geq \min\{n,2\delta\}$.\\

When  $G\backslash C$ has a simple structure, namely is edgeless, Voss and 
Zuluaga [6] obtained the following.\\

\noindent\textbf{Theorem C} [6]. Let $G$ be a 3-connected graph. Then either 
$c\geq 3\delta-3$ or each longest cycle in $G$ is a dominating cycle.\\

The first lower bound involving connectivity $\kappa$ as a parameter has been appeared in 1981, 
by the author [3].\\

\noindent\textbf{Theorem D} [3]. Let $G$ be a 2-connected graph. Then 
$c\geq \min\{n,3\delta-\kappa\}$.\\

Further, the first two bounds involving  $\overline{p}$ and $\overline{c}$ has been appeared in 1998 and 2000, respectively, again by the author [4],[5].\\  

\noindent\textbf{Theorem E} [4]. Let $G$ be a graph and $C$ a longest cycle in $G$. 
Then $|C|\geq(\overline{p}+2)(\delta-\overline{p})$.\\

\noindent\textbf{Theorem F} [5]. Let $G$ be a graph and $C$ a longest cycle in $G$. 
Then $|C|\geq(\overline{c}+1)(\delta-\overline{c}+1)$.\\

As a defect, the bound in Theorem D decreases as $\kappa$ 
increases. The bounds in Theorems E and F have the same defect for $\overline{p}\geq(\delta-2)/2$
 and $\overline{c}\geq\delta/2$, respectively.

In this paper we present a lower bound for the circumference involving $\delta$, 
$\kappa$ and $\overline{c}$ and increasing with $\delta$, $\kappa$ and $\overline{c}$. \\

\noindent{\bf Theorem 1.}\quad Let $G$ be a graph and $C$ a longest cycle in $G$. Then
$\left| C\right| \geq \frac{(\overline{c}+1)\kappa}{\overline{c}+\kappa+1}\left( \delta+2\right)$ if 
$\overline{c}\geq \kappa$ and $|C|\geq
\frac{(\overline{c}+1)\overline{c}}{2\overline{c}+1}(\delta+2)$ if $\overline{c}\leq\kappa-1$.
.\\

The result is sharp, as can be
seen from the following family of graphs. Take $\kappa+1$ disjoint copies of the
complete graph $K_{\delta -\kappa+1}$ and join each vertex in their union to
every vertex of a disjoint complete graph $K_\kappa$. This graph $(\kappa+1)K_{\delta
-\kappa+1}+K_\kappa$ is clearly not hamiltonian. Moreover, $c=\kappa(\delta -\kappa+2)$ and $%
\overline{c}=\delta -\kappa+1$, implying that $c=\frac{(\overline{c}+1)\kappa}{\overline{c}+\kappa+1}\left( \delta +2\right)
.$\\

In view of Theorem 1, we belive the following is also true in terms of $\overline{p}$.\\

\noindent\textbf{Conjecture 1}. Let $G$ be a graph and $C$ a longest cycle in $G$.
Then $|C|\geq\frac{(\overline{p}+2)\kappa}{\overline{p}+\kappa+2}(\delta +2)$ if $\overline{p}\geq \kappa-1$, and
$|C|\geq\frac{(\overline{p}+2)\overline{p}}{2\overline{p}+2}(\delta +2)$ if $\overline{p}\le \kappa-2$.
\\

The next section is devoted to standard terminology. In section 3 we
introduce some special definitions and convenient notations, where the
notion of $HC-$extensions plays a central role in the sequel. In section 4
we investigate the main properties of $HC-$extensions and in the last
section we prove our main result.

\section{Terminology}

An $(x,y)$-path is a path with endvertices $x$ and $y$. Given an $(x,y)$%
-path $L$ of $G$, we denote by $\overrightarrow{L}$ the path $L$ with an
orientation from $x$ to $y$. If $u,v\in V(L)$, then $u\overrightarrow{L}v$
denotes the consequtive vertices on $L$ from $u$ to $v$ in the direction
specified by $\overrightarrow{L}$. The same vertices, in reverse order, are
given by $v\overleftarrow{L}u$. For $\overrightarrow{L}=x\overrightarrow{L}y$
and $u\in V(L)$, let $u^{+}(\overrightarrow{L})$ (or just $u^{+}$) denotes
the successor of $u$ $(u\not =y)$ on $\overrightarrow{L}$ and $u^{-}$
denotes its predecessor $(u\not =x)$. If $A\subseteq V(L)$, then $%
A^{+}=\left\{ v^{+}\mid v\in A\backslash \{y\}\right\} $ and $A^{-}=\left\{ v^{-}\mid v\in
A\backslash \{x\}\right\} $. If $Q$ is a cycle in $G$ and $A\subseteq V(Q)$, then $%
\overrightarrow{Q}$, $A^{+}$ and $A^{-}$ are analogously defined. For $v\in
V(Q)$, $v\overrightarrow{Q}v$ will be interpreted as a vertex v. For $v\in V$, 
put $N(v)=\left\{ u\in V\mid uv\in E\right\} $ and $d(v)=\left|N(v)\right| $.

\section{Special Definitions}

We begin introducing some special definitions and convenient notations. For
the remainder of this section let a longest cycle $C$ in a graph $G$ and a
longest cycle $H=u_1\ldots u_hu_1$ in $G\backslash C$ with $h=\overline{c}$ be fixed.\\

\noindent{\bf Definition 3.1.}\quad $HC-$extension; $T(u_i);\stackrel{o}{u}%
;\hat u$.

Let $T(u_1),\ldots ,T(u_h)$ are vertex-disjoint $(u_i,\hat u_i)$-paths in $%
G\backslash C$ $(i=1,...,h)$. The union $T=\bigcup_{i=1}^hT(u_i)$
is called $HC-$extension if $N(\hat u_i)\subseteq V(T)\bigcup V(C)$ for each $
i\in \{1,...,h\}$. An HC-extension $T$ is called maximal if it is chosen
so as to maximize $\left| \left\{ u\in V(H)\left| u\neq \hat u\right.
\right\} \right| $. If $u\neq \hat u$ for some $u\in V(H)$, then we use $%
\stackrel{o}{u}$ to denote $u^{+}(\overrightarrow{T}(u))$.\\

\noindent{\bf Definition 3.2.\quad }$(A,B)-$path.

Let $A,B\subset V$ and $A\bigcap B=\emptyset $. Let $E$ is a path in $G$ with all
its inner vertices in $V\backslash (A\bigcup B)$. Then $E$ is called an $(A,B)$-path if $E$ starts
at any vertex in $A$ and terminates at any vertex in $B$. For subgraphs $H_1$
and $H_2$ of $G$, an $(H_1,H_2)$-path is analogously defined.\\

\noindent{\bf Definition 3.3.}\quad $\Theta (\overrightarrow{P},V_{neut},V_{fin})=(P_0,...,P_\pi)$; 
$P_i=y_i\overrightarrow{P_i}z_i$ $(i=0,\ldots ,\pi )$.

Let $V^{^{\prime }}\subset V$. A path with endvertices in $V\backslash V^{^{\prime }}$
and all internal vertices in $V^{^{\prime }}$ is called a $V^{^{\prime }}-$%
path. Let $\overrightarrow{P}=v_0v_1\ldots v_n$ be a path in $G$ of length $%
n\geq 1$ and let $V_{neut}$, $V_{fin}$ be vertex-disjoint subsets in $V\backslash V(%
\overrightarrow{P})$. We define $\Theta (\overrightarrow{P}%
,V_{neut},V_{fin}) $ as a sequence of paths $P_0,\ldots ,P_\pi $ as follows:
For $i=0$, put $P_0=\overrightarrow{y_0z_0}$ and $X=V(v_0\overrightarrow{P}%
z_0)$, where $y_0=v_0$ and $z_0=v_1$. Now let $P_{i-1}=y_{i-1}%
\overrightarrow{P}_{i-1}z_{i-1}$ and $X_{i-1}$ are defined for some integer $%
i\geq 1$. In order to define $P_i$ and $X_i$ we distinguish three cases.

$(i)$ If every $V_{neut}-$path, starting in $X_{i-1}-z_{i-1}$, terminates in 
$X_{i-1}$, then $X_i=\emptyset $ and $P_\pi =P_{i-1}$ (so $P_i$ is undefined).

$(ii)$ If there is a $V_{neut}-$path $P^{\prime }=v^{\prime }\overrightarrow{%
P^{\prime }}v^{\prime \prime }$ with $v^{\prime }\in X_{i-1}-z_{i-1}$ and $%
v^{\prime \prime }\in V_{fin}$, then $X_i=\emptyset $ and $P_\pi =P_i=y_i%
\overrightarrow{P^{\prime }}z_i$, where $y_i=v^{\prime }$ and $z_i=v^{\prime
\prime }$.

$(iii)$ There is a $V_{neut}-$path $P^{\prime \prime }=w^{\prime }%
\overrightarrow{P^{\prime \prime }}w^{\prime \prime }$ with $w^{\prime }\in
X_{i-1}-z_{i-1}$ and $w^{\prime \prime }\in V(z_{i-1}^{+}\overrightarrow{P}%
v_n)$ but there is no $V_{neut}-$path satisfying (ii).

Choose $P^{\prime \prime }$ so as to maximize $\mid v_0\overrightarrow{P}%
w^{\prime \prime }\mid $. Then putting 
$$
P_i=y_i\overrightarrow{P^{\prime \prime }}z_i,~X_i=V(v_0\overrightarrow{P}%
z_i), 
$$
where $y_i=w^{\prime }$ and $z_i=w^{\prime \prime }$, we complete the
definition of $P_i$ and $X_i$. Since $X_o\subset X_1\subset \cdots $, there
must be some integer $j$ $(j\geq 1)$ with $P_j=P_\pi $, which, in fact,
completes the definition of $\Theta (\overrightarrow{P},V_{neut},V_{fin})$.\\

\noindent{\bf Definition 3.4.} $\Phi _u;\varphi _u;\Psi _u;\psi _u$.

Let $T$ be a maximal HC-extension. For each $u\in V(H)$, put 
$$
\Phi _u=N(\hat u)\bigcap V(T),\quad \varphi _u=\left| \Phi _u\right| , 
$$
$$
\Psi _u=N(\hat u)\bigcap V(C),\quad \psi _u=\left| \Psi _u\right| . 
$$

\noindent{\bf Definition 3.5.} $U_0;\overline{U}_0;U_1;U_2;U_{*}$.

For $T$ a maximal $HC-$extension, put 
$$
U_0=\{u\in V(H)\left| u=\hat u\right. \},\enskip \overline{U}_0=V(H)\backslash U_0,\enskip
U_1=\{u\in \overline{U}_0\left| \Phi _u\not \subseteq V(T(u))\right. \}. 
$$
Let $u\in V(H)\backslash(U_0\bigcup U_1)$ and $\Theta (\overleftarrow{T}%
(u),V_{neut},V_{fin})=(P_0,\ldots ,P_\pi )$, where 
$$
V_{neut}=V\backslash(V(T)\bigcup V(C)),\quad V_{fin}=V(T)\backslash V(T(u)). 
$$
A vertex $u$ is called special if $P_\pi $ starts and terminates in $%
V(T(u))$. The set of all nonspecial vertices in $V(H)\backslash(U_0\bigcup U_1)$ is
denoted by $U_2$ and the set of all special vertices by $U_{*}$.\\

\noindent{\bf Definition 3.6.}$\quad B_u;B_u^{*};b_u;b_u^{*}$.

Let $T$ be a maximal HC-extension. For each $u\in V(H)$, put $B_u=\{v\in
U_0\mid v\stackrel{o}{u}\in E\}$. Clearly $B_u=\emptyset $ if $u\in U_0$.
Furthermore, for each $u\in U_0$, put $B_u^{*}=\{v\in V(H)\mid u\stackrel{o}{%
v}\in E$. Clearly $B_u^{*}\subseteq \overline{U}_0$. Let $b_u=\left|
B_u\right| $ and $b_u^{*}=\left| B_u^{*}\right| $.\\

\noindent{\bf Definition 3.7.\quad }$A_u(v);\rho _u(v);\overline{\rho }_u(v);\Lambda
_u;\Lambda _u(v,w)$.

Let $T$ be a maximal HC-extension. For each $u,v\in V(H)$, put 
$$A_u(v)=(\Phi_u\bigcup B_u)\bigcap V(T(v)).$$ 
Let  $\rho _u(v)$ denote the vertex in  $%
A_u(v)$ maximizing  $\mid v\overrightarrow{T}(v)\rho _u(v)\mid $. In
particular, $\rho _u(u)=\hat u^{-}$. Put $\overline{\rho }_u(v)=\hat u$ if $%
\rho _u(v)\in \Phi _u$ and $\overline{\rho }_u(v)=\stackrel{o}{u}$ if $\rho
_u(v)\in B_u$. Clearly $\overline{\rho }_u(u)=\hat u$. Put $\Lambda
_u=\left\{ v\in V(H)\left| A_u(v)\not =\emptyset \right. \right\} $. For
each $v,w\in \Lambda _u\quad (v\not =w)$, put 
$$
\Lambda _u(v,w)=vT(v)\rho _u(v)\overline{\rho }_u(v)T(u)\overline{\rho }%
_u(w)\rho _u(w)T(w)w. 
$$

\noindent{\bf Definition 3.8.\quad }$\varphi _u^{\prime };\gamma _u;\beta _u;\mu (T).$

For $T$ a maximal $HC-$extension, put\\

$\varphi _u^{\prime }=\left\{ 
\begin{array}{lll}
\varphi _u & \mbox{if} & \mbox{ }u\in V(H)\backslash U_{*}, \\ 0 & \mbox{if} & \mbox{ }%
u\in U_{*}, 
\end{array}
\right. \qquad \gamma _u=\left\{ 
\begin{array}{lll}
\varphi _u^{^{^{\prime }}}+b_u & \mbox{if} & \mbox{ }u\in \overline{U}_0, \\ 
\varphi _u^{^{\prime }}-b_u^{*} & \mbox{if} & \mbox{ }u\in U_0, 
\end{array}
\right. $
\\ \\
$$
\beta _u=\frac{(\gamma _u+\gamma _{u^{+}})}2\quad (u\in V(H)),\qquad \mu
(T)=\frac{1}{h}\sum_{u\in V(H)}\beta _u. 
$$

\noindent{\bf Definition 3.9.\quad }$T-$transformation; $%
T_{tr}(E_1,...,E_n);T_{tr}(v_1,...,v_n).$

Let $T$ be a maximal $HC-$extension and let $E_1,...,E_n$ are
vertex-disjoint $(H,C)-$paths with $E_i=v_i\overrightarrow{E}_iw_i\quad
(i=1,...,n)$. Assume that the union of $E_1,...,E_n$ intersect $T(z)$ for
some $z\in V(H)\backslash\{v_1,...,v_n\}$. Clearly $z\in \overline{U}_0$. walking
along $T(z)$ from $z$ to $\stackrel{\wedge }{z}$ we stop at the first vertex 
$w\in \bigcup_{i=1}^nV(E_i)$. Assume w.l.o.g. that $w\in V(E_1)$. Replacing
the segment $v_1E_1w$ of a path $E_1$ by $zT(z)w$ we get a new path $E_1^0$
instead of $E_1$. If the union of $E_1^0,E_2,...,E_n$ intersect $%
T(z^{^{\prime }})$ for some $z^{^{\prime }}\in V(H)\backslash\{z,v_2,...,v_n\}$, then
continue this procedure. In a finite number of steps we obtain 
$$
\mid \{v\in V(H): \big(\bigcup^n_{i=1} 
V(E_i^{^{\prime }})\big)\bigcap V(T(v))\not =\emptyset \}\mid =n 
$$
for some vertex-disjoint $(H,C)$-paths $E_1^{^{\prime }},...,E_n^{^{\prime
}}.$ Let $E_i^{^{\prime }}=v_i^{^{\prime }}E_i^{^{\prime }}w_i\quad
(i=1,...,n).$ By writing 
$$
T_{tr}\left( E_1,...,E_n\right) =(E_1^{^{\prime }},...,E_n^{^{\prime
}}),\quad T_{tr}\left( v_1^{^{}},...,v_n^{}\right) =\left( v_1^{^{\prime
}},...,v_n^{^{\prime }}\right) , 
$$
we say that $E_1^{^{\prime }},...,E_n^{^{\prime }}$ is a $T-$ transformation
of $E_1,...,E_n$. By the definition,

$$
v_i^{^{\prime }}\in \left\{ v_i\right\} \bigcup \overline{U}_0\quad \left(
i=1,...,n\right) ,\quad T_{tr}\left( w_1,...,w_n\right) =\left(
w_1,...,w_n\right) . 
$$

\noindent{\bf Definition 3.10.}\quad $O(x,y);O_x(x,y);O(y,\stackrel{o}{x});O(x,\stackrel{o}{%
x});O_y(x,y);O(x,\stackrel{o}{y});O(y,\stackrel{o}{y}).$

Let $T$ be a maximal $HC-$extension. For each pair of distinct vertices $%
x,y\in V(H),$ put%
$$
V_1=\bigcup_{v\not \in \left\{ x,y\right\} }V\left( T\left(
v\right) \right) \bigcup \left\{ x,y\right\} ,\quad V_2=V_1\bigcup \{%
\stackrel{o}{x}\}. 
$$

Let $O(x,y)$ (resp. $O_x(x,y),O(y,\stackrel{o}{x}),O(x,\stackrel{o}{x}%
)) $ be the longest $(x,y)-$path (resp. $(x,y)$-path, $(y,\stackrel{o}{x})-$%
path, $(x,\stackrel{o}{x})-$path) in $\langle V_1\rangle$
(resp. $\langle V_2\rangle$,$\langle V_2\rangle$,%
$\langle V_2\rangle)$. The paths $O_y(x,y),$ $O(x,\stackrel{o}{y%
}),$ and $O(y,\stackrel{o}{y})$ are analogously defined.\\

\noindent{\bf Definition 3.11.}\quad $\Omega \left( x,y\right) ;\Omega \left(
x,y,E,F\right) ;\Omega \left( v,w,x,y,E,F\right) .$

Let $T$ be a maximal $HC-$extension and let $E,F$ be a pair of vertex
disjoint $T-$trans\-for\-med $(H,C)-$paths with $E=xEv$ and $F=yFw.$ If $\left|
T\left( x\right) \right| -1\not =1$, then we denote $\Omega _x\left(
x,y,E,F\right) =O(x,y).$ Otherwise,%
$$
\Omega _x\left( x,y,E,F\right) =\left\{ 
\begin{array}{lll}
O_x(x,y) & \mbox{if} & \stackrel{o}{x}\not \in V(E)\bigcup V(F), \\ O(%
\stackrel{o}{x},y) & \mbox{if} & \stackrel{o}{x}\in V(E), \\ O(\stackrel{o}{x%
},x) & \mbox{if} & \stackrel{o}{x}\in V(F). 
\end{array}
\right. 
$$

Define $\Omega _y\left( x,y,E,F\right) $ by the same way and denote by $\Omega
\left( x,y,E,F\right) $ the longest path among $O(x,y),\Omega _x\left(
x,y,E,F\right) $ and $\Omega _y\left( x,y,E,F\right) .$ Let $\Omega \left(
x,y\right) $ be the shortest path $\Omega \left( x,y,E,F\right) $ for fixed $%
x,y$ and all possible $E,F.$ By definition 3.9, $vE\mu \Omega \left(
x,y,E,F\right) \nu Fw$ is a simple path for appropriate $\mu ,\nu \in \{x,y,%
\stackrel{o}{x},\stackrel{o}{y}\}$ denoted by $\Omega \left(
v,w,x,y,E,F\right) .$\\

\noindent{\bf Definition 3.12.\quad }$(v,L)\in \Delta .$

Let $L$ be a path in $G$ with $L=v_1...v_{2t-1}\quad (t\geq 1)$ and let $%
v\in V\backslash V(L).$ We will write $(v,L)\in \Delta $ if $vv_{2i-1}\in E\quad
(i=1,...,t).$ For $w\in V(L)$, we will write $(w,L)\in \Delta $ if $wu\in
E$ for each $u\in V(L)\backslash \{w\}.$\\

\noindent{\bf Remarks. }If no ambiguity can arise, any notation of the type $R_{u_i}$
in definitions 2.4 and 2.6-2.8, having index $u_i$ (say $\Phi _{u_i}$), we
abbreviate $R_{u_i}=R_i.$

\section{Preliminaries}

Throughout this section, let $G$ be a graph, $C$ be a longest cycle $G$ and $%
H=u_1...u_hu_1$ a longest cycle in $G\backslash C$ with a maximal $HC-$extension $T.$\\

\noindent{\bf Lemma 1.} Let $G$ be a graph.

{\bf (a1)}\quad For $E,F$ a pair of vertex-disjoint $(H,C)-$paths in $G$ with $%
E=xEv$ and $F=yFw$, if $T_{tr}\left( E,F\right) =\left( E^{\prime
},F^{\prime }\right) $ and $T_{tr}\left( x,y\right) =\left( x^{\prime
},y^{\prime }\right) $, then 
$$
\mid v\overrightarrow{C}w\mid -1\geq \left| \Omega (v,w,x^{\prime
},y^{\prime },E^{\prime },F^{\prime })\right| -1\geq a(x^{\prime
})+a(y^{\prime })+\left| \Omega (x^{\prime },y^{\prime })\right| -1, 
$$
where $a(z)=1$ if $z\not \in U_{*}$ and $a(z)=\varphi _z+1$ if $z\in U_{*}$
for each $z\in \left\{ x^{\prime },y^{\prime }\right\} .$

{\bf (a2)} Let $u\in V(H)$ and $\Theta (\overleftarrow{T}%
(u),V_{neut},V_{fin})=(P_0,...,P_\pi ),$ where $P_i=y_i\overrightarrow{P}%
_iz_i\quad (i=0,...,\pi )$ and 
$$
V_{neut}=V\backslash(V(T)\bigcup V(C)),\quad V_{fin}=V(T)\backslash V(T(u)). 
$$

If $u\in U_2$, then there is an $(u,z_\pi )-$path $L$ of length at least $%
\varphi _u+1$ with $V(L)\subseteq V(T(u))\bigcup V^{*}$, where $%
V^{*}=\bigcup_{i=0}^\pi V(P_i).$ If $u\in U_{*}$, then for each vertex 
$$
z\in (V(\hat u\overleftarrow{T}(u)z_\pi )\bigcup V^{*})\backslash \{z_\pi\} 
$$
there is an $\left( u,z\right) -$path $L$ of length at least $\varphi _u+1$
with $V\left( L\right) \subseteq V\left( T\left( u\right) \right) \bigcup
V^{*}.$\\

\noindent{\bf Lemma 2.} For each $u\in V(H)$,

$
\begin{array}{ll}
{\bf (b1)}  \quad  \mbox{if} \enskip u\in 
\overline{U}_0 \enskip \mbox{and} \enskip \hat u\not =\stackrel{o}{u}, \enskip \mbox{then} \enskip \Phi _u\bigcap
B_u=\emptyset .\\
{\bf (b2)} \quad  \sum_{u\not\in U_0}b_u=\sum_{u\in U_0}b_u^*,  
\enskip \sum_{u\in V(H)}\gamma _u=\sum_{u\in V(H)}\varphi _u^{\prime },\\
\left| \Phi _u\bigcup B_u\right| =\sum_{v\in V(H)}\left| A_u\left( v\right) \right|. 
\end{array}
$\\

\noindent{\bf Lemma 3.} Let $G$ be a graph, $C$ be a longest cycle in $G$, $Q$ be a path in $G\backslash C$
and let $P_i=v_i\overrightarrow{P_i}w_i\quad (i=0,...,q)$ are vertex-disjoint
paths in $G\backslash C$ having only $v_0,...,v_q$ in common with $Q.$ Then%
$$
c\geq \sum^q_{i=0}\left| Z_i\right| +\mid \bigcup^q_{i=0}Z_i\mid , 
$$
where $Z_i=N(w_i)\bigcap V(C)\quad (i=0,...,q).$\\

\noindent{\bf Lemma 4.} For each $u\in V(H),$

$
\begin{array}{llll}
{\bf (d1)} \quad \mbox{if} \enskip \left| T\left( u\right) \right| -1\geq 2, \enskip\mbox{then} \enskip 
h\geq 2\gamma _u. \\ 
{\bf (d2)}  \quad \mbox{if}\enskip\left| T\left( u\right) \right| -1=1,\enskip\mbox{then} \enskip  h\geq
2\varphi _u^{\prime }\geq \gamma _u+1. \\ 
{\bf (d3)}  \quad h\geq \gamma _u+1.     
\end{array}
$\\

\noindent{\bf Lemma 5.} Let $\Lambda _u\subseteq V(x\overrightarrow{H}y)$ for some $%
u,x,y\in V(H).$

$
\begin{array}{l}
{\bf (e1)\quad }\mbox{if}\enskip\left| T\left( u\right) \right| -1\geq 2,\enskip\mbox{then}\enskip
 \mid x\overrightarrow{H}y\mid -1\geq \gamma _u. \\ 
{\bf (e2)\quad }\mbox{if}\enskip \left| T\left( u\right) \right| -1=1, \enskip\mbox{then}\enskip 
\mid x\overrightarrow{H}y\mid -1\geq \gamma _u-1. \\ 
{\bf (e3)\quad }\mbox{if}\enskip\left| T\left( u\right) \right| -1=1 \enskip\mbox{and}\enskip\mid x\overrightarrow{H}y\mid -1=\gamma
_u-1,\enskip\mbox{then}\\ 
 (\hat u,x\overrightarrow{H}y)\in \Delta,\enskip
B_u=\Lambda _u  \backslash \{u\}\subseteq U_0\enskip\mbox{and}\enskip\gamma _u-1=2(\varphi _u-1). 
\end{array}
$\\

\noindent{\bf Lemma 6.} For each $u\in U_1\bigcup U_2$, let $x_1\overrightarrow{H}y_1$
and $x_2\overrightarrow{H}y_2$ be vertex-disjoint segments in $H$ with $%
\left\{ x_1,x_2,y_1,y_2\right\} \subseteq \Lambda _u\subseteq V(x_1%
\overrightarrow{H}y_1)\bigcup V(x_2\overrightarrow{H}y_2)$ and let $v\in
\left\{ x_2,y_2\right\} .$

$
\begin{array}{ll}
{\bf (f1)} & \mbox{If }B_u\bigcup \left\{ u\right\} \subseteq V(x_1%
\overrightarrow{H}y_1)\mbox{ and }\Lambda _u\backslash\left( B_u\bigcup \left\{
u\right\} \right) \subseteq V(x_2\overrightarrow{H}y_2),\mbox{ then} 
\end{array}
$%
$$
\mid x_1\overrightarrow{H}y_1\mid -1+\mid x_2\overrightarrow{H}y_2\mid
-1+\left| A_u\left( u\right) \right| +\left| A_u\left( v\right) \right| \geq
\gamma _u-1. 
$$
Otherwise,%
$$
\mid x_1\overrightarrow{H}y_1\mid -1+\mid x_2\overrightarrow{H}y_2\mid
-1+\left| A_u\left( u\right) \right| +\left| A_u\left( v\right) \right| \geq
\gamma _u-2+\left| A_u\left( u\right) \right| \geq \gamma _u-1. 
$$

$
\begin{array}{ll}
{\bf (f2)} & \mbox{If }\mid x_1\overrightarrow{H}y_1\mid -1+\mid x_2%
\overrightarrow{H}y_2\mid -1+\left| A_u\left( u\right) \right| +\left|
A_u\left( v\right) \right| =\gamma _u-1,\mbox{ then} 
\end{array}
$%
$$
(\hat u,x_i\overrightarrow{H}y_i)\in \Delta \mbox{ }(i=1,2),\quad
B_u=\Lambda _u\backslash \{u\}\subseteq U_0,\quad \gamma _u-1=2(\varphi _u-1). 
$$

\noindent{\bf Lemma 7.} Let $x,y$ be a pair of distinct vertices in $H$. For each $%
u\in V(H),$

${\bf (g1)\quad }\mbox{if} \enskip u\in U_{*}, \enskip\mbox{then} \enskip  \left| O(x,y)\right|
-1\geq \gamma _u+1$,

${\bf (g2)\quad }\mbox{if}\enskip\left| T\left( u\right) \right| -1\geq 2,\enskip\mbox{then}\enskip
 \left| O(x,y)\right| -1\geq \gamma _u$,

${\bf (g3)\quad }\mbox{if}\enskip\left| T\left( u\right) \right| -1=1, \enskip\mbox{then}\enskip
\left| O_u(x,y)\right| -1\geq \gamma _u-1$,

${\bf (g4)}$\quad Let $\left| T\left( u\right) \right| -1=1$ and $\left|
O_u(x,y)\right| -1=\gamma _u-1.$ If either $\Lambda _u\subseteq V(x%
\overrightarrow{H}y)$ and ($\stackrel{o}{u},H)\not \in \Delta $ or $\Lambda
_u\subseteq V(y\overrightarrow{H}x)$ and
$(\stackrel{o}{u},H)\not \in \Delta $ (say $\Lambda _u\subseteq V(x%
\overrightarrow{H}y)$ and $(\stackrel{o}{u},H)\not \in \Delta )$, then

$\quad {\bf (g4.1)\quad }(\stackrel{o}{u},x\overrightarrow{H}y)\in \Delta $,

$\quad {\bf (g4.2)\quad }B_u=\Lambda _u\backslash \{u\}\subseteq U_0\enskip\mbox{and}\enskip\left|
O_u(x,y)\right| -1=\mid x\overrightarrow{H}y-1\mid =\gamma _u-1=2(\varphi
_u-1),$

$\quad {\bf (g4.3)\quad }\mbox{if} \enskip z\in V(x\overrightarrow{H}y)\backslash\Lambda
_u, \enskip\mbox{then \enskip either}\enskip z\in U_{*}\enskip \mbox{or}\enskip 
z\in U_0 \enskip \mbox{and} \\ \Lambda_z\subseteq 
\Lambda _u\bigcup \left\{ z\right\}, \enskip \gamma _z\leq \varphi _u=\left( \gamma _u+1\right) /2,$

$\quad {\bf (g4.4)\quad } \mbox{if} \enskip z\in V(x\overrightarrow{H}y)\backslash\left\{ x,y\right\}, \enskip\mbox{then}\enskip
\Lambda _z\subseteq V(x\overrightarrow{H}y).$

Otherwise,

$\quad {\bf (g4.5)\quad }(\stackrel{o}{u},H)\in \Delta ,$

$\quad {\bf (g4.6)\quad }B_u=\Lambda _u\backslash \{u\}\subseteq U_0 \enskip\mbox{and}\enskip \left|
O_u(x,y)\right| -1=h-2=\gamma _u-1=2(\varphi _u-1),$

$\quad {\bf (g4.7)\quad } \mbox{if}\enskip z\in V\left( H\right) \backslash\Lambda _u,\enskip\mbox{then either}\enskip
 z\in U_{*} \enskip \mbox{or}\enskip z\in U_0,\enskip\Lambda _z\subseteq \Lambda _u\bigcup
\left\{ z\right\}\\ \mbox{and} \enskip \gamma _z\leq \varphi _u=(\gamma _u+1)/2=h/2,$

${\bf (g5)\quad }  \mbox{if} \enskip \left| T\left( x\right) \right| -1=1, \enskip\mbox{then}\enskip  \min
\{\left| O(\stackrel{o}{x},x)\right| -1,\left| O(\stackrel{o}{x},y)\right|
-1\}\geq \gamma _x.$

${\bf (g6)\quad } \mbox{if}\enskip u\in \left\{ x^{+},x^{-},y^{+},y^{-}\right\},
\enskip\mbox{then}\enskip\left| O(x,y)\right| -1\geq \gamma _u,$

${\bf (g7)\quad }$If $u\in \left\{ x^{+},x^{-},y^{+},y^{-}\right\} \quad ($%
say $u=x^{+})$ and $\left| O(x,y)\right| -1=\gamma _u$, then $\mid T\left(
u\right) \mid -1\leq 1$
and $(\hat u,v\overrightarrow{H}y)\in \Delta $ for some $v\in \Lambda _u$
with $\Lambda _u\subseteq V(v\overrightarrow{H}y),$

${\bf (g8)\quad }$If $\left| T\left( x\right) \right| -1=1$ and $\left|
O_x(x,y)\right| -1=\left| O_x(x,w)\right| -1=\gamma _x-1$ for some

$w\in V\left( H\right) \backslash\left\{ x,y\right\} $, then for each $z\in \left\{
x^{+},x^{-}\right\} ,$
$$
\min \{\left| O_x(x,y)\right| -1,\left| O_x(x,w)\right| -1\}\geq \gamma _z+1%
\mbox{.} 
$$

\noindent{\bf Lemma 8. }Let $x,y$ be a pair of distinct vertices in $H$ and let 
$$
\begin{array}{c}
a=\min \{\left| O_x(x,y)\right| ,\mid O( 
\stackrel{o}{x},y)\mid ,\mid O(\stackrel{o}{x},x)\mid\} -1, \\ b=\min \{\left|
O_y(x,y)\right| ,\mid O(\stackrel{o}{y},x)\mid ,\mid O(\stackrel{o}{y}%
,y)\mid \}-1. 
\end{array}
$$
Then $\left| \Omega (x,y)\right| -1\geq \max \{\left| O(x,y)\right| -1,a,b\}.$\\

\noindent{\bf Lemma 9. }Let $x,y$ be a pair of distinct vertices in $H.$ Then

${\bf (i1)\quad } \mbox{if} \enskip \left\{ u_i,u_{i+1}\right\} \bigcap \left\{ x,y\right\}
=\emptyset ~(i\in \{1,...,h\}), \enskip\mbox{then}\enskip \ \left| \Omega
(x,y)\right| -1\geq (\gamma _i+\gamma _{i+1})/2=\beta _i,$

${\bf (i2)\quad } \mbox{if}\enskip \left| T\left( x\right) \right| -1\geq 2 \enskip\mbox{and}\enskip z\in \left\{
x^{+},x^{-}\right\}, \enskip\mbox{then}\enskip \ \left| \Omega (x,y)\right| -1\geq
(\gamma _x+\gamma _z)/2$,

${\bf (i3)\quad } \mbox{if}\enskip x\in U_{*} \enskip\mbox{and}\enskip z\in \left\{ x^{+},x^{-}\right\}, 
\enskip\mbox{then}\enskip \left| \Omega (x,y)\right| -1\geq (\gamma _x+\gamma _z+1)/2$,

${\bf (i4)\quad }$If $\left| T\left( x\right) \right| -1=1$, then for each $%
w\in V\left( H\right) \backslash\left\{ x,y\right\} $ and $z\in \left\{
x^{+},x^{-}\right\} ,$%
$$
\max \{\left| \Omega (x,y)\right| -1,\left| \Omega \left( x,w\right) \right|
-1\}\geq (\gamma _x+\gamma _z)/2, 
$$

${\bf (i5)\quad } \mbox{if} \enskip z\in \left\{ x^{+},x^{-}\right\} \enskip\mbox{and}\enskip w\in V\left( H\right)
\backslash \{z\}, \enskip\mbox{then}\enskip \max \{\left| \Omega (x,y)\right| -1,\\ \left| \Omega \left(
z,w\right) \right| -1\}\geq (\gamma _x+\gamma _z)/2,$

${\bf (i6)\quad }$If $x\in \overline{U}_0$ and $h\not =4$, then $\left|
\Omega (x,y)\right| -1\geq (\gamma _x+\gamma _z)/2$ for some $z\in \left\{
x^{+},x^{-}\right\}$,

${\bf (i7)\quad } \mbox{if}\enskip x,y\in \overline{U}_0, \enskip\mbox{then}\enskip \ \left| \Omega \left(
x,y\right) \right| -1\geq \max_i\beta _i,$

${\bf (i8)\quad } \mbox{if} \enskip \mid x\overrightarrow{H}y\mid -1=1, \enskip\mbox{then}\enskip  \ \left|
\Omega \left( x,y\right) \right| -1\geq \max_i\beta _i,$

${\bf (i9)\quad } \mbox{if} \enskip \mid x\overrightarrow{H}y\mid -1=2 \enskip\mbox{and}\enskip ~h\not =4, \enskip\mbox{then}\enskip
\left| \Omega (x,y)\right| -1\geq (\gamma _x+\gamma _{x^{+}})/2.\ $

\section{Proofs}

{\bf Proof of lemma 1.} {\bf (a1)} Following definition 3.11, we distinguish
three cases.

{\bf Case 1.} $x,y\not \in U_{*}.$

Clearly, \\

$
\begin{array}{ll}
\mid v\overrightarrow{C}w\mid -1\geq \mid \Omega (v,w,x,y,E,F)\mid -1\\
\geq 2+\mid \Omega (x,y,E,F)\mid -1\geq \mid \Omega (x,y)\mid +1. 
\end{array}
$\\

{\bf Case 2.} $x,y\in U_{*}.$

If $\left| T\left( x\right) \right| -1=1$, then $\stackrel{o}{x}\not \in
V\left( \Omega \left( x,y\right) \right) ,$ since otherwise the segment of $%
\Omega \left( x,y\right) $ between $\stackrel{o}{x}$ and $y$ , contradict
the fact that $x\in U_{*}.$ Therefore, $\Omega \left( x,y,E,F\right)
=O\left( x,y\right) .$ On the other hand, $\Omega _x\left( x,y,E,F\right)
=O\left( x,y\right) $ if $\left| T\left( x\right) \right| -1\geq 2.$ Also, by
the symmetric arguments, $\Omega _y\left( x,y,E,F\right) =O\left( x,y\right)
.$ Thus $\Omega \left( x,y,E,F\right) =\Omega \left( x,y\right) $ and\\

$
\begin{array}{ll}
\mid v 
\overrightarrow{C}w\mid -1\geq \left| \Omega \left( v,w,x,y,E,F\right)
\right| -1 \\ \geq \left( \left| E\right| -1\right) +\left( \left| F\right|
-1\right) +\left| \Omega \left( x,y\right) \right| -1\geq \varphi _x+\varphi
_y+\left| \Omega \left( x,y\right) \right| +1. 
\end{array}
$

{\bf Case 3.} Either $x\not \in U_{*},y\in U_{*}$ or $x\in U_{*},y\not \in
U_{*}.$

Apply the arguments in case 1 and case 2. \quad $\Delta $\\

${\bf (a2)}~$Suppose first that $u\in U_2.$ By definition 3.3, $z_1\in
V\left( T\left( u\right) \right) $ and $z_\pi \in V_{fin}.$ Let $z_\pi \in
V\left( T\left( w\right) \right) $ for some $w\in V\left( H\right) \backslash \{\}u.$
Choose $z_{01}\in V(u\overrightarrow{T}\left( u\right) y_2^{-})$ such that $%
z_{01}\hat u\in E$ and $\mid z_{01}\overrightarrow{T}\left( u\right) y_2\mid 
$ is minimum. Then we get the desired result putting together the following
paths%
$$
P_2,...,P_\pi,\enskip\hat uz_{01},\enskip\hat u\overleftarrow{T}\left( u\right) y_2,\enskip z_{01}%
\overleftarrow{T}\left( u\right) y_3,\enskip z_{\pi -1}\overleftarrow{T}\left(
u\right) u, \enskip z_i\overleftarrow{T}\left( u\right) y_{i+2}, 
$$
where $i=2,...,\pi-2$. A similar proof holds for $u\in U_{*}$. \qquad $\Delta $\\

\noindent{\bf Proof of lemma 2. (b1)}. {\bf Case 1.} $u\in U_1.$

Suppose, to the contrary, that $\Phi _u\bigcap B_u\not =\emptyset $ and let 
$z\in $ $\Phi _u\bigcap B_u.$ Then, by definitions 3.4 and 3.1, the
collection%
$$
\{T\left( u_1\right) ,...,T\left( u_h\right) ,u\stackrel{o}{u},z\hat
u\}\backslash\left\{ T\left( u\right) ,T\left( z\right) \right\} 
$$
generates another $HC-$extension, contradicting the maximality of $T$.

{\bf Case 2.} $u\in U_2\bigcup U_{*}.$

By definition 3.5, $\Phi _u\subseteq V\left( T\left( u\right) \right) $ and
the result follows. \quad $\Delta $

${\bf (b2)}~$Immediately from definitions 3.6-3.8. \quad $\Delta $\\

\noindent{\bf Proof of lemma 3.} Assume first that $v_i=w_i~\left( i=0,...,q\right) $%
. The result is immediate if $\bigcup_{i=0}^qZ_i=\emptyset $. Let $%
\bigcup_{i=0}^qZ_i\not =\emptyset $ and let $\xi _1,...,\xi _m~\left( m\geq
1\right) $ be the elements of $\bigcup_{i=0}^qZ_i$ occuring on $%
\overrightarrow{C}$ in a consequtive order. Set 
$$
F_i=N\left( \xi _i\right) \bigcap \left\{ w_0,...,w_q\right\} \quad \left(
i=1,...,m\right) . 
$$
Suppose that $m=1$. If $\left| F_1\right| =1$, then $q=0$ and $%
Z_0=Z_q=\left\{ \xi _1\right\} $ implying that 
$$
c\geq 2=\sum^q_{i=0}\left| Z_i\right| +\mid 
\bigcup^q_{i=0}Z_i\mid . 
$$

If $\left| F_1\right| \geq 2$, then choosing $u,v\in F_1~(u\not =v)$ such
that $\mid u\overrightarrow{Q}v\mid $ is maximum, we get
$$
c\geq \mid \xi _1u\overrightarrow{Q}v\xi _1\mid \geq \sum^q_{i=0}\left| Z_i\right| +1
=\sum^q_{i=0}\left| Z_i\right| +\mid \bigcup^q_{i=0}Z_i\mid . 
$$

Thus, we may assume $m\geq 2$. It means, in particular, that $c\geq 3$. For $%
i=1,...,m$ , put $f\left( \xi _i\right) =$ $\mid \xi _i\overrightarrow{C}\xi
_{i+1}\mid -1~\left( \mbox{indices mod }m\right) $. It is easy to see that 
\begin{equation}
\label{1}
\begin{array}{lll}
c=\sum^m_{i=1}f\left( \xi _i\right) , & \sum^m_{i=1}\left| F_i\right| =\sum^q_{i=0}\left| Z_i\right| ,
 & m=\mid \bigcup^q_{i=0}Z_i\mid . 
\end{array}
\end{equation}
For every $i\in \{1,...,m\}$, choose $x_i,y_i\in F_i\bigcup F_{i+1}$ such
that $\mid x_i\overrightarrow{Q}y_i\mid $ is maximum (indices mod $m$).\\

{\bf Claim 3.1}\quad $f\left( \xi _i\right) \geq \left( \left| F_i\right|
+\left| F_{i+1}\right| +2\right) /2\quad \left( i=1,...,m\right) .$

{\bf Proof of Claim 3.1.} We distingwish two cases.

{\bf Case 1.}\quad Either $x_i\in F_i,~y_i\in
F_{i+1}$ or $x_i\in F_{i+1},~y_i\in F_i.$

If $x_i\in F_i,~y_i\in F_{i+1}$, then $f\left( \xi _i\right) \geq $ $\mid \xi
_ix_i\overrightarrow{Q}y_i\xi _{i+1}\mid -1$ and hence%
$$
f\left( \xi _i\right) \geq \max \left( \left| F_i\right| ,\left|
F_{i+1}\right| \right) +1\geq \left( \left| F_i\right| +\left|
F_{i+1}\right| +2\right) /2. 
$$
Otherwise, the result holds from $f\left( \xi _i\right) \geq \mid \xi _iy_i%
\overleftarrow{Q}x_i\xi _{i+1}\mid -1$ in the same way.

{\bf Case 2. }Either $x_i,~y_i\in F_i$ or $x_i,~y_i\in F_{i+1}.$

First, suppose $x_i,~y_i\in F_i.$ We can assume also $x_i,~y_i\not \in
F_{i+1} $, since otherwise we could argue as in case 1. Choose $%
x_i^{^{\prime }},~y_i^{^{\prime }}\in F_{i+1}$ such that $\mid x_i^{^{\prime
}}\overrightarrow{Q}y_i^{^{\prime }}\mid $ is maximum. If $\mid x_i%
\overrightarrow{Q}x_i^{^{\prime }}\mid -1\geq \left( \left| F_i\right|
-\left| F_{i+1}\right| \right) /2$, then%
$$
f(\xi _i)\geq \mid \xi _ix_i\overrightarrow{Q}y_i^{^{\prime }}\xi _{i+1}\mid
-1\geq \left( \left| F_i\right| -\left| F_{i+1}\right| \right) /2+\left|
F_{i+1}\right| +1= \left( \left| F_i\right| +\left| F_{i+1}\right|
+2\right) /2. 
$$
Otherwise,\\

$
\begin{array}{ll}
f(\xi _i)\geq \mid \xi _iy_i 
\overleftarrow{Q}x_i^{^{\prime }}\xi _{i+1}\mid -1=\mid x_i^{^{\prime }}%
\overrightarrow{Q}y_i\mid +1=\mid x_i\overrightarrow{Q}y_i\mid -\mid x_i%
\overrightarrow{Q}x_i^{^{\prime }}\mid +2\geq \\ \geq \left| F_i\right|
-\left( \left| F_i\right| -\left| F_{i+1}\right| +1\right) /2+2> \left(
\left| F_i\right| +\left| F_{i+1}\right| +2\right) /2. 
\end{array}
$\\

By symmetry, the case $x_i,~y_i\in F_{i+1}$ requires the same arguments.$ \quad
\Delta $\\

By claim 3.1,%
$$
\sum^m_{i=1}f\left( \xi _i\right) \geq \sum^m_{i=1}\left( \left| F_i\right| +\left| F_{i+1}\right|
+2\right) /2=\sum^m_{i=1}\left| F_i\right| +m, 
$$
which by (\ref{1}) gives the desired result. Finally, if $v_i\not =w_i$ for
some $i\in \{0,...,q\}$, then we can argue exactly as in case $%
v_i=w_i\quad (i=0,...,q)$.  \qquad $\Delta $\\

\noindent{\bf Proof of lemma 4.} {\bf (d1)}. {\bf Case1.}\quad $u\in U_1.$

Let $\xi _1,...,\xi _f$ be the elements of $\Lambda _u$  
occuring on $H$ in a consequtive order with $u=\xi _1$. For
each integer $i~\left( 1\leq i\leq f\right) $, let%
$$
\begin{array}{lll}
M_i=\xi _i\overrightarrow{H}\xi _{i+1}, & \omega _i=\left| A_u\left( \xi
_i\right) \right| +\left| A_u\left( \xi _{i+1}\right) \right| & \left( \mbox{%
indices mod }f\right) 
\end{array}
. 
$$
Since $H$ is extreme, 
\begin{equation}
\label{2}\left| M_i\right| \geq \left| \Lambda _u\left( \xi _i,\xi
_{i+1}\right) \right| \quad \left( i=1,...,f\right) . 
\end{equation}
Let $\xi _r\overrightarrow{H}\xi _s$ be the longest segment on $H$ such that
$$
\begin{array}{ll}
\xi _1\in V(\xi _r\overrightarrow{H}\xi _s), & \left\{ \xi _r,\xi
_{r+1},...,\xi _s\right\} \subseteq B_u\bigcup \left\{ u\right\} . 
\end{array}
$$
Put%
$$
\begin{array}{l}
\Omega ^{+}=\left\{ M_i\in \left\{ M_2,...,M_{f-1}\right\} \left| 
\overline{\rho }_u\left( \xi _i\right) \not =\overline{\rho }_u\left( \xi
_{i+1}\right) \right. \right\} , \\ \Omega ^{-}=\left\{ M_i\in \left\{
M_1,M_f\right\} \left| 
\overline{\rho }_u\left( \xi _i\right) \not =\overline{\rho }_u\left( \xi
_{i+1}\right) \right. \right\} , \\ \Omega ^0=\left\{ M_1,...,M_f\right\}
\backslash(\Omega ^{+}\bigcup \Omega ^{-}). 
\end{array}
$$
Observe that $\left| \Omega ^{-}\right| \leq 2$ and $\left| M_i\right|
-1\geq \left| \Lambda _u\left( \xi _i,\xi _{i+1}\right) \right| -1$ for each 
$i\in \{1,...,f\}.$ Then clearly 
\begin{equation}
\label{3}
\begin{array}{lll}
\mbox{if} \enskip M_i\in \Omega ^{+}, \enskip\mbox{then}\enskip \left| M_i\right| -1\geq \omega
_i+\left| A_u\left( u\right) \right| -1, 
\end{array}
\end{equation}
\begin{equation}
\label{4}
\begin{array}{lll}
\mbox{if} \enskip M_i\in \Omega ^{-}, \enskip\mbox{then}\enskip \left| M_i\right| -1\geq \omega
_i-\left| A_u\left( u\right) \right| +1, 
\end{array}
\end{equation}
\begin{equation}
\label{5}
\begin{array}[t]{lll}
M_i\in \Omega ^0 & \Longrightarrow & \left| M_i\right| -1\geq \omega _i. 
\end{array}
\end{equation}

{\bf Claim 4.1.} If $\left| \Omega ^{-}\right| =0$ then $\left| M_i\right|
-1\geq \omega _i~\left( i=1,...,f\right) $.

{\bf Proof of Claim 4.1.} Immediate from (\ref{3}), (\ref{4}) and (\ref{5}).$\quad
\Delta $\\

{\bf Claim 4.2.}

${\bf (k1)\quad }$If $\left| \Omega ^{-}\right| =1$,
say $\Omega ^{-}=\left\{ M_1\right\} ,$ then $M_s\in \Omega ^{+}.$

${\bf (k2)\quad }$If $\Omega ^{-}=\left\{ M_1\right\} $ and $\Omega
^{+}=\left\{ M_s\right\}$, then%
$$
\begin{array}{ll}
B_u\bigcup \left\{ u\right\} \subseteq V(\xi _r\overrightarrow{H}\xi _s), & 
\Lambda _u\backslash\left( B_u\bigcup \left\{ u\right\} \right) \subseteq V(\xi
_{s+1}\overrightarrow{H}\xi _{r-1}). 
\end{array}
$$

{\bf Proof of Claim 4.2.}\quad${\bf (k1)\ }$Let $\Omega ^{-}=\left\{
M_1\right\} .$ By the definition, $\left\{ \xi _2,...,\xi _s\right\}
\subseteq B_u$ and $\xi _{s+1}\in \Lambda _u\backslash\left( B_u\bigcup \left\{
u\right\} \right) ,$ implying that $M_s\in \Omega ^{+}$.\quad$\Delta $\\

${\bf (k2)\quad }$If follows that $\left\{ \xi _{s+1},...,\xi _f\right\}
\subseteq \Lambda _u\backslash\left( B_u\bigcup \left\{ u\right\} \right) .$ On the
other hand (by the definition), $\left\{ \xi _1,...,\xi _s\right\} \subseteq
B_u\bigcup \left\{ u\right\} $ and the proof is complete.\quad$\Delta $\\

{\bf Claim 4.3.}\quad ${\bf (l1)\quad }$If $\left| \Omega ^{-}\right| =2,$
i.e.$\Omega ^{-}=\left\{ M_1,M_f\right\}$, then $M_s,M_{r-1}\in \Omega
^{+}. $

${\bf (l2)~}$If $\Omega ^{-}=\left\{ M_1,M_f\right\} $ and $\Omega
^{+}=\left\{ M_s,M_{r-1}\right\}$, then $\xi _1,\xi _r,\xi _s$ are pairwise
different and 
$$
\begin{array}{cc}
B_u\bigcup \left\{ u\right\} \subseteq V(\xi _r\overrightarrow{H}\xi _s), & 
\Lambda _u\backslash\left( B_u\bigcup \left\{ u\right\} \right) \subseteq V(\xi
_{s+1}\overrightarrow{H}\xi _{r-1}). 
\end{array}
$$

$\Pr ${\bf oof of claim 4.3.} ${\bf (l1)\quad }$By the definition, $\left\{
\xi _2,\xi _f,\xi _s,\xi _r\right\} \subseteq B_u$ and $\xi _{s+1},\xi
_{r-1}\in \Lambda _u\backslash\left( B_u\bigcup \left\{ u\right\} \right) $, which
implies $M_s,M_{r-1}\in \Omega ^{+}$. \quad$\Delta $\\

${\bf (l2)~}$It follows that $\left\{ M_{s+1},...,M_{r-2}\right\} \bigcap
\Omega ^{+}=\emptyset $ and hence%
$$
\left\{ \xi _{s+1},...,\xi _{r-1}\right\} \subseteq \Lambda _u\backslash\left(
B_u\bigcup \left\{ u\right\} \right) . 
$$
On the other hand (by the definition) $\left\{ \xi _r,...,\xi _s\right\}
\subseteq B_u\bigcup \left\{ u\right\} $ , which completes the proof of
claim 4.3. \quad $\Delta $\\

The following three statements can be proved easely basing on (\ref{3}), (\ref{4}%
), (\ref{5}) and claims 4.1, 4.2 and 4.3.\\

{\bf Claim 4.4.}\quad $\sum_{i=1}^f\left( \left| M_i\right| -1\right) \geq
\sum_{i=1}^f\omega _i.$\\

{\bf Claim 4.5.}\quad \mbox{If} \enskip $t\in \left\{ 1,...,f\right\}, \enskip\mbox{then}\enskip  \sum_{i%
\not =t}\left( \left| M_i\right| -1\right) \geq \sum_{i\not =t}\omega
_i-\left| A_u\left( u\right) \right| +1.$\\

{\bf Claim 4.6.}\quad \mbox{If} \enskip $g,t\in \left\{ 1,...,f\right\} ~\left( g\not
=t\right), \enskip\mbox{then}\enskip  \sum_{i\not \in \left\{ g,t\right\} }\left( \left|
M_i\right| -1\right) \geq\\ \sum_{i\not \in \left\{ g,t\right\} }\omega
_i-2\left| A_u\left( u\right) \right| +2.$\\

Using $(b2)$ and claim 4.4, we get\\

$
\begin{array}{ll}
h= 
\sum^f_{i=1}\left( \left| M_i\right| -1\right) \geq 
\sum^f_{i=1}\omega _i=\sum^f_{i=1}\left( \left| A_u\left( \xi _i\right) \right|
 +\left| A_u\left( \xi_{i+1}\right) \right| \right) \\ =2\sum^f_{i=1}\left| 
A_u\left( \xi _i\right) \right| =2\left| \Phi _u\bigcup B_u\right| . 
\end{array}
$\\

By $(b1)$, $\left| \Phi _u\bigcup B_u\right| =\varphi _u+b_u=\gamma _u,$
implying that $h\geq 2\gamma _u.$

{\bf Case 2.}\quad $u\in U_2.$

Let $\Theta (\overleftarrow{T}\left( u\right) ,V_{neut},V_{fin})=\left(
P_0,...,P_\pi \right) $, where $P_i=y_i\overrightarrow{P_i}z_i~\left(
i=0,...,\pi \right) .$ By $(a2)$, there is an $\left( u,z_\pi \right) -$path 
$L$ of length at least $\varphi _u+1$ with $V\left( L\right) \subseteq
V\left( T\left( u\right) \right) \bigcup V^{*}.$ Let $z_\pi \in V(T(w))$
for some $w\in V(H)$. By denoting $B_u\bigcup \left\{ u,w\right\}\\
 =\left\{ \xi_1,...,\xi _f\right\}$, we can argue exactly as in case 1.

{\bf Case 3.}\quad $u\in U_{*}.$

Clearly $h\geq 2\left( b_u+1\right) =2(\varphi _u^{^{\prime
}}+b_u+1)>2\gamma _u. \quad \Delta $\\

${\bf (d2)~}$Since $\left| T\left( u\right) \right| -1=1$, we have $u\in
U_1\bigcup U_{*}.$If $u\in U_{*}$, then $b_u=0$ and $h\geq 2=2(\varphi
_u^{^{\prime }}+b_u+1)=2\left( \gamma _u+1\right) \geq \gamma _u+1$. Let $%
u\in U_1$. Define 
\begin{equation}
\label{6}\xi _i,\omega _i,M_i,\xi _r\overrightarrow{H}\xi _s,\Omega
^{+},\Omega ^{-},\Omega ^0 
\end{equation}
as in proof of $(d1){\bf .}$ It is easy to see that $\Omega ^{+}=\Omega
^{-}=\emptyset .$ By claim 4.1, $\sum_{i=1}^f\left( \left| M_i\right|
-1\right) \geq \sum_{i=1}^f\omega _i$ and as in proof of $(d1){\bf ,}$ $%
h\geq 2\left| \Phi _u\bigcup B_u\right| =2\varphi _u.$ Noting that $\varphi
_u\geq b_u+\left| \left\{ u\right\} \right| =b_u+1$, we obtain $h\geq
\varphi _u+b_u+1=\gamma _u+1. \quad\Delta $\\

${\bf (d3)}$ It is easely checked that $h\geq \gamma _u+1$ if $u\in U_0$. If 
$u\in \overline{U}_0$, then by $(d1)$ and $(d2)$, $h\geq \min \left(
2\gamma _u,\gamma _u+1\right) \geq \gamma _u+1.$         \quad  $\Delta $\\

\noindent{\bf Proof of lemma 5.} Assume w.l.o.g. that $x,y\in \Lambda _u$.\\ 

${\bf (e1)\quad }${\bf Case 1\quad }$u\in U_1.$

Following (\ref{6}) we let, in addition, $y\overrightarrow{H}x=M_t$ for some 
$t\ \left( 1\leq t\leq f\right) $. By claim 4.5 we can distingush the
following two cases:

{\bf Case 1.1.}$\quad \mid x\overrightarrow{H}y\mid -1\geq \sum_{i\not
=t}\omega _i.$

By $(b1){\bf ,}\left| \Phi _u\bigcup B_u\right| =\left| \Phi _u\right|
+\left| B_u\right| =\gamma _u$, and using $(b2),$
$$
\begin{array}{c}
\mid x 
\overrightarrow{H}y\mid -1\geq \sum_{i\not =t}\omega _i=\sum^f_{i=1}\omega _i-\omega _t\\ =2
\sum ^f_{i=1}\left| A_u\left( \xi _i\right) \right| -\left| A_u\left( \xi
_t\right) \right| -\left| A_u\left( \xi _{t+1}\right) \right| \\ =\sum^f_{i=1}\left| 
A_u\left( \xi _i\right) \right| +\sum_{i\not\in\{t,t+1\}}\left|A_u(\xi_i)\right|\\
 \geq \sum^f_{i=1}\left| A_u\left(
\xi _i\right) \right| =\left| \Phi _u\bigcup B_u\right| =\gamma _u. 
\end{array}
$$

{\bf Case 1.2.}\quad $\sum_{i\not =t}\omega _i-\left| A_u\left( u\right)
\right| +1\leq \left| x\overrightarrow{H}y-1\right| <\sum_{i\not =t}\omega _i.$

If $\Omega ^{-}=\emptyset $, then by claim 4.1., $\mid x\overrightarrow{H}%
y\mid -1\geq \sum_{i\not =t}\omega _i,$ a contradition. Let $\Omega ^{-}\not
=\emptyset $.

{\bf Case 1.2.1.}\quad $\left| \Omega ^{-}\right| =1.$

Assume w.l.o.g $\Omega ^{-}=\left\{ M_1\right\} .$ By claim 4.2, $M_s\in
\Omega ^{+}.$ If $\left| \Omega ^{+}\right| \geq 2$, then by (\ref{3}), (\ref
{4}) and (\ref{5}), $\mid x\overrightarrow{H}y\mid -1\geq \sum_{i\not
=t}\omega _i,$ a contradiction. Thus we can assume $\Omega ^{+}=\left\{
M_s\right\} .$ If $M_t\not =M_s$, then again $\mid x\overrightarrow{H}y\mid
-1\geq \sum_{i\not =t}\omega _i,$ a contradiction. Finally, if $M_t=M_s$,
then $A_u\left( \xi _t\right) $, $A_u\left( \xi _{t+1}\right) $ and $%
A_u\left( u\right) $ are pairwise different and hence

$$
\begin{array}{c}
\mid x 
\overrightarrow{H}y\mid -1\geq 2\sum^f_{i=1}\left|
A_u\left( \xi _i\right) \right| -\left| A_u\left( \xi _t\right) \right|
-\left| A_u\left( \xi _{t+1}\right) \right| -\left| A_u\left( u\right)
\right| +1= \\ \sum^f_{i=1}\left| A_u\left( \xi
_i\right) \right| +\sum_{i\not \in \left\{ 1,t,t+1\right\} }\left| A_u\left(
\xi _i\right) \right| +1\geq \sum^f_{i=1}\left|
A_u\left( \xi _i\right) \right| +\left( f-3\right) +1 \\ \geq \sum^f_{i=1}
\left| A_u\left( \xi _i\right) \right| =\left| \Phi
_u\bigcup B_u\right| =\gamma _u. 
\end{array}
$$

{\bf Case 1.2.2.}\quad $\left| \Omega ^{-}\right| =2.$

By claim 4.3, $M_s,M_{r-1}\in \Omega ^{+}$. If $\left| \Omega ^{+}\right|
\geq 3$, then by (\ref{3}), (\ref{4}) and (\ref{5}), $\mid x\overrightarrow{H}%
y\mid -1\geq \sum_{i\not =t}\omega _i,$ a contradiction. Let $\Omega ^{+}$ $%
=\left\{ M_s,M_{r-1}\right\} $. If $M_t\not \in \Omega ^{+}$, then again $%
\mid x\overrightarrow{H}y\mid -1\geq \sum_{i\not =t}\omega _i,$ a
contradiction. Finally, if $M_t\in \Omega ^{+},$ say $M_t=M_s$, then $%
A_u\left( \xi _t\right) $, $A_u\left( \xi _{t+1}\right) ,$ $A_u\left(
u\right) $ are pairwise different and we can argue exactly as in case
1.2.1.

{\bf Case 2.}\quad $u\in U_2\bigcup U_{*}.$

Apply the arguments used in the proof of $(d1){\bf \quad }\left( \mbox{case
2 and case 3}\right) .$\\

${\bf (e2)\quad }$Clearly $u\in U_1$. Following (\ref{6}) we see that $%
\Omega ^{+}=\Omega ^{-}=\emptyset $. By claim 4.1, $\left| M_i\right| -1\geq
\omega _i~\left( i=1,...,f\right) $. Recalling that $f\geq b_u+1$,%
$$
\begin{array}{c}
\mid x 
\overrightarrow{H}y\mid -1=\sum_{i\not =t}\left( \left|
M_i\right| -1\right) \geq \sum_{i\not =t}\omega _i=\sum^f_{i=1}
\left| A_u\left( \xi _i\right) \right| +\sum_{i\not
\in \left\{ t,t+1\right\} }\left| A_u\left( \xi _i\right) \right| \\ \geq
\left| \Phi _u\bigcup B_u\right| +f-2=\varphi _u+f-2=\varphi
_u+b_u-1=\gamma _u-1. 
\end{array}
$$

${\bf (e3)\quad }$It was shown in $(e2)$ that $\mid x\overrightarrow{H}y\mid
-1\geq \varphi _u+f-2\geq $ $\gamma _u-1$. Since $\mid x\overrightarrow{H}%
y\mid -1=\gamma _u-1$, we have $\mid x\overrightarrow{H}y\mid -1=\varphi
_u+f-2=\gamma _u-1$. This implies $\left| B_u\right| =b_u=f-1$ and therefore 
$B_u=\Lambda _u\backslash \{u\}\subseteq U_0$. But then $\varphi _u=b_u+1$ and $\mid x%
\overrightarrow{H}y\mid -1=\gamma _u-1=2\varphi _u-2$ implying that ($\hat
u,x\overrightarrow{H}y)\in \Delta $. \qquad $\Delta $\\

{\bf Proof of lemma 6.} ${\bf (f1)}$ {\bf Case 1.\quad }$u\in U_1.$

By symmetry, we can assume $v=x_2$. Following (\ref{6}) we let, in addition, 
$M_g=y_1\overrightarrow{H}x_2$ and $M_t=y_2\overrightarrow{H}x_1$ for some
integers $g,t\in \left\{ 1,...,f\right\} $. This means that%
$$
y_1=\xi _g,~x_2=\xi _{g+1},~y_2=\xi _t,~x_1=\xi _{t+1},~v=x_2=\xi
_{g+1},~A_u\left( \xi _{g+1}\right) =A_u\left( v\right) . 
$$
Putting $\beta =\mid x_1\overrightarrow{H}y_1\mid -1+\mid x_2\overrightarrow{%
H}y_2\mid -1$ and using claim 4.6, we can distingush the following four
cases.

{\bf Case 1.1.}\quad $\beta \geq \sum_{i\not \in \left\{ g,t\right\} }\omega
_i+\mid A_u\left( u\right) \mid -1.$

Clearly
\begin{equation}
\label{7}
\begin{array}{c}
\beta \geq 2 
\sum^f_{i=1}\left| A_u\left( \xi _i\right) \right|
-\left| A_u\left( \xi _g\right) \right| -\left| A_u\left( \xi _{g+1}\right)
\right| -\left| A_u\left( \xi _t\right) \right| -\left| A_u\left( \xi
_{t+1}\right) \right| \\ +\left| A_u\left( u\right) \right| -1. 
\end{array}
\end{equation}
Observe that $\left| A_u\left( u\right) \right| \geq 1$ and $f\geq b_u+1$.
If $x_1\not =y_1$, then $A_u\left( \xi _g\right)$, $\ A_u\left( \xi _t\right)$
and $\ ~A_u\left( \xi _{t+1}\right) $ are pairwise different and by (\ref{7}),
$$
\begin{array}{lll}
\beta +\left| A_u\left( v\right) \right| =\beta +\left| A_u\left( \xi
_{g+1}\right) \right| \geq 
\sum^f_{i=1}\left| A_u\left( \xi _i\right) \right| +\sum
_{i\not \in \left\{ g,t,t+1\right\} }\left| A_u\left( \xi
_i\right) \right|  \\\geq \sum^f_{i=1}\left|
A_u\left( \xi _i\right) \right| +f-3\geq \left| \Phi _u\bigcup B_u\right|
+f-3 \\\geq \left| \Phi _u\right| +f-3\geq \varphi _u+b_u-2\geq \gamma
_u-1-\left| A_u\left( u\right) \right| . 
\end{array}
$$

Otherwise $\left( x_1=y_1\right) ,~A_u\left( \xi _{t+1}\right) =A_u\left(
u\right) $ , and by (\ref{7}),
$$
\begin{array}{lll}
\beta +\left| A_u\left( v\right) \right| =\beta +\left| A_u\left( \xi
_{g+1}\right) \right| \geq 2 
\sum^f_{i=1}\left| A_u\left( \xi _i\right) \right|
-\left| A_u\left( \xi _g\right) \right| -\left| A_u\left( \xi _t\right)
\right| -1 \\ \geq 
\sum^f_{i=1}\left| A_u\left( \xi _i\right) \right|
+\sum_{i\not \in \left\{ g,t\right\} }\left| A_u\left( \xi _i\right) \right|
-1\geq \left| \Phi _u\bigcup B_u\right| +f-3\\ \geq \varphi _u+b_u-2\geq  
\gamma _u-1-\left| A_u\left( u\right) \right| . 
\end{array}
$$

{\bf Case 1.2.}\quad $\sum_{i\not \in \left\{ g,t\right\} }\omega _i\leq
\beta <\sum_{i\not \in \left\{ g,t\right\} }\omega _i+\left| A_u\left(
u\right) \right| -1.$

Clearly
$$
\beta +\left| A_u\left( v\right) \right| =\beta +\left| A_u\left( \xi
_{g+1}\right) \right|
$$
$$
 \geq 2\sum^f_{i=1}\left|
A_u\left( \xi _i\right) \right| -\left| A_u\left( \xi _g\right) \right|
-\left| A_u\left( \xi _t\right) \right| -\left| A_u\left( \xi _{t+1}\right)
\right| . 
$$

If $x_1\not =y_1$, then we obtain the desired result as in case 1.1. Let $%
x_1=y_1$, i.e. $M_g=M_1,$ $M_t=M_f$ and $\xi _{t+1}=\xi _g=\xi _1=u$. If $%
\Omega ^{+}\not =\emptyset $, then $\beta \geq \sum_{i\not \in \left\{
1,t\right\} }\omega _i+\left| A_u\left( u\right) \right| -1$, a
contradiction. Let $\Omega ^{+}=\emptyset $. This implies $M_s=M_1$ and $%
M_{r-1}=M_f$, and we deduce that%
$$
B_u\bigcup \left\{ u\right\} =\left\{ u\right\} =V(x_1\overrightarrow{H}%
y_1),\quad \Lambda _u\backslash \{u\}\subseteq V(x_2\overrightarrow{H}y_2). 
$$
Recalling that $f\geq b_u+1$, we get
$$
\begin{array}{c}
\beta +\left| A_u\left( u\right) \right| +\left| A_u\left( v\right) \right|
=\beta +\left| A_u\left( u\right) \right| +\left| A_u\left( \xi
_{g+1}\right) \right|\\ \geq 2 
\sum^t_{i=1}\left| A_u\left( \xi _i\right) \right|
-\left| A_u\left( \xi _t\right) \right|  -\left| A_u\left( \xi
_{t+1}\right) \right|\\ \geq \left| \Phi _u\bigcup B_u\right| +f-2\geq
\varphi _u+b_u-1\geq \gamma _u-1. 
\end{array}
$$

{\bf Case 1.3.}\quad $\sum_{i\not \in \left\{ g,t\right\} }\omega _i-\mid
A_u(u)\mid +1\leq \beta <\sum_{i\not \in \left\{ g,t\right\} }\omega _i.$\\

{\bf Case 1.3.1.}\quad $\xi _1\not \in \left\{ x_1,y_1\right\} .$

It follows that $A_u\left( u\right) ,A_u\left( \xi _g\right) ,A_u\left( \xi
_t\right) $ and $A_u\left( \xi _{t+1}\right) $ are pairwise different. Since 
$f\geq b_u+1,$ we have\\

$
\begin{array}{lll}
\beta +\left| A_u\left( v\right) \right| =\beta +\left| A_u\left( \xi
_{g+1}\right) \right|\\ \geq 2 
\sum^f_{i=1}\left| A_u\left( \xi _i\right) \right|
-\left| A_u\left( u\right) \right| -\left| A_u\left( \xi _g\right) \right|
-\left| A_u\left( \xi _t\right) \right|  -\left| A_u\left( \xi
_{t+1}\right) \right| +1 \\ \geq \sum^f_{i=1}\left|
A_u\left( \xi _i\right) \right| +f-3\geq \varphi _u+b_u-2\geq \gamma
_u-1-\left| A_u\left( u\right) \right| . 
\end{array}
$\\

{\bf Case 1.3.2.}\quad $\xi _1\in \left\{ x_1,y_1\right\} .$

Assume w.l.o.g. that $\xi _1=x_1$, i.e. $M_t=M_f$. If $M_1\not \in \Omega
^{-}$, then $\beta \geq \sum_{i\not \in \left\{ g,t\right\} }\omega _i$, a
contradiction. Let $M_1\in \Omega ^{-}$. This implies $\xi _2\in B_u$ and $%
M_s\in \Omega ^{+}$. If $M_s\not =M_g$, then $\beta \geq \sum_{i\not \in
\left\{ g,t\right\} }\omega _i$, a contradiction. So, assume $M_s=M_g$.
Analogously, $M_{r-1}=M_f$. If $M_j\in \Omega ^{+}$ for some $j\in \left\{
1,...,f\right\} \backslash\left\{ g\right\} $, then again $\beta \geq \sum_{i\not \in
\left\{ g,t\right\} }\omega _i$, a contradiction. Let%
$$
i\in \left\{ 1,...,f\right\} \backslash\left\{ g\right\} ~\Longrightarrow ~M_i\in
\Omega ^0\bigcup \Omega ^{-}. 
$$

It follows that $B_u\bigcup \left\{ u\right\} \subseteq V(x_1%
\overrightarrow{H}y_1)$ and $\Lambda _u\backslash(B_u\bigcup \left\{ u\right\}
)\subseteq V(x_2\overrightarrow{H}y_2)$. Furthermore, noting that $A_u\left(
\xi _g\right) ,A_u\left( \xi _t\right) ,A_u\left( \xi _{t+1}\right) $ are
pairwise different and $f\geq b_u+1$, we get
$$
\begin{array}{lll}
\beta +\left| A_u\left( u\right) \right| +\left| A_u\left( v\right) \right|
=\beta +\left| A_u\left( u\right) \right| +\left| A_u\left( \xi
_{g+1}\right) \right| \\ \geq 2 
\sum^f_{i=1}\left| A_u\left( \xi _i\right) \right|
-  \left| A_u\left( \xi _g\right) \right| -\left| A_u\left( \xi _t\right)
\right| -\left| A_u\left( \xi _{t+1}\right) \right| +1 \\  \geq 
\sum^f_{i=1}\left| A_u\left( \xi _i\right) \right|
+f-2\geq  \varphi _u+b_u-1\geq \gamma _u-1. 
\end{array}
$$

{\bf Case 1.4.}\quad $\sum_{i\not \in \left\{ g,t\right\} }\omega _i-2\left|
A_u\left( u\right) \right| +2\leq \beta <\sum_{i\not \in \left\{ g,t\right\}
}\omega _i-\left| A_u\left( u\right) \right| +1.$

If $\left| \Omega ^{-}\right| \leq 1$, then clearly $\beta \geq \sum_{i\not
\in \left\{ g,t\right\} }\omega _i-\left| A_u\left( u\right) \right| +1$, a
contradiction. Let $\left| \Omega ^{-}\right| =2$. This implies $%
~M_1,~M_f\in \Omega ^{-}$ and $M_s,M_{r-1}\in \Omega ^{+}$. If $\left|
\Omega ^{+}\right| \geq 3$, then again $\beta \geq \sum_{i\not \in \left\{
g,t\right\} }\omega _i-\left| A_u\left( u\right) \right| +1$, a
contradiction. Let $\left| \Omega ^{+}\right| =2$, i.e. $\Omega ^{+}=\left\{
M_s,M_{r-1}\right\} $. By claim 4.3,%
$$
B_u\bigcup \left\{ u\right\} \subseteq V(x_1\overrightarrow{H}y_1),\quad
\Lambda _u\backslash(B_u\bigcup \left\{ u\right\} )\subseteq V(x_2\overrightarrow{H}%
y_2). 
$$
Recalling that $f\geq b_u+1$, we get
$$
\begin{array}{lll}
\beta +\left| A_u\left( u\right) \right| +\left| A_u\left( v\right) \right|
=\beta +\left| A_u\left( u\right) \right| +\left| A_u\left( \xi
_{g+1}\right) \right| \\ 
\geq 2 
\sum^f_{i=1} \left| A_u\left( \xi _i\right) \right|
-\left| A_u\left( u\right) \right| -\left| A_u\left( \xi _g\right) \right|
-\left| A_u\left( \xi _t\right) \right| -\left| A_u\left( \xi _{t+1}\right)
\right| +2 \\ \geq \sum ^f_{i=1}\left| A_u\left( \xi
_i\right) \right| + \sum _{i\not\in\{1,g,t,t+1\}}
\left| A_u\left( \xi _i\right) \right| +2 \\ \geq \left| \Phi _u\bigcup
B_u\right| +\left( f-4\right) +2\geq \gamma _u-1. 
\end{array}
$$

{\bf Case 2.}\quad $u\in U_2.$

Apply the arguments used in the proof of $(d1)$ (see case 2 and case 3).\\

${\bf (f2)\quad }${\bf Case 1.\quad }$u\in U_1.$

As shown in the proof of $(f1){\bf ,}$%
$$
\beta +\left| A_u\left( u\right) \right| +\left| A_u\left( v\right) \right|
\geq \left| \Phi _u\bigcup B_u\right| +f-2\geq \left| \Phi _u\right|
+f-2\geq \varphi _u+b_u-1. 
$$
Since $\beta +\left| A_u\left( u\right) \right| +\left| A_u\left( v\right)
\right| =\varphi _u+b_u-1=\gamma _u-1,$ we have equations%
$$
\beta +\left| A_u\left( u\right) \right| +\left| A_u\left( v\right) \right|
=\left| \Phi _u\bigcup B_u\right| +f-2=\left| \Phi _u\right| +f-2=\varphi
_u+b_u-1 
$$
implying that $f=b_u+1$. If $\left| T\left( u\right) \right| -1\geq 2$, then $%
\Lambda _u\backslash(B_u\bigcup \left\{ u\right\} )\not =\emptyset $ and hence $%
f\geq \left| B_u\right| +\left| \left\{ u\right\} \right| +1=b_u+2$, a
contradiction. Otherwise $(\left| T\left( u\right) \right| -1=1)$,%
$$
\begin{array}{c}
\begin{array}{ll}
\Lambda _u=B_u\bigcup \left\{ u\right\} , & \left| A_u\left( u\right)
\right| =\left| A_u\left( v\right) \right| =1, 
\end{array}
\\ 
\begin{array}{ll}
\varphi _u=b_u+1, & \beta =\varphi _u+b_u-3=2\varphi _u-4 
\end{array}
\end{array}
$$
and we deduce that $(\hat u,x_i\overrightarrow{H}y_i)\in \Delta ~\left(
i=1,2\right) $ and $\gamma _u-1=2\left( \varphi _u-1\right) $.

{\bf Case 2.}\quad $u\in U_2.$

Apply the arguments used in the proof of $(d1){\bf}$. \qquad $\Delta $\\

{\bf Proof of lemma 7.} ${\bf (g1)~}$Clearly $h\geq 2\left( b_u+1\right)
=2(\varphi _u^{^{\prime }}+b_u+1)=2\left( \gamma _u+1\right) $ and
therefore, $\left| O(x,y)\right| -1\geq h/2\geq \gamma _u+1$.\quad$\Delta$\\

${\bf (g2)\quad }$ By (d1), $h\geq 2\gamma _u$ which implies $\left|
O(x,y)\right| -1\geq h/2\geq \gamma _u.\quad \Delta $\\

${\bf (g3),(g5)\quad }$If $\left| \Lambda _u\right| =1$, then $\gamma _u=1$
and there is nothing to prove. Let $\left| \Lambda _u\right| \geq 2$, i.e. $%
u\in U_1$.

{\bf Case 1.}\quad $u\not \in \left\{ x,y\right\} .$

Assume w.l.o.g. that $u\in V(x^{+}\overrightarrow{H}y^{-}).$ We can assume
also that $\Lambda _u\not \subseteq V(x\overrightarrow{H}y)$, since
otherwise the result holds by $(e2)$. Let $x_1\overrightarrow{H}y_1$ be the
longest segment in $x\overrightarrow{H}y^{-}$ with $x_1,y_1\in \Lambda _u$
and $x_2\overrightarrow{H}y_2$ be the longest segmnet in $y\overrightarrow{H}%
x^{-}$ with $x_2,y_2\in \Lambda _u$. Putting $\beta =\mid x_1\overrightarrow{%
H}y_1\mid -1+$ $\mid x_2\overrightarrow{H}y_2\mid -1$ we see (by lemma 6)
that%
$$
\beta \geq \gamma _u-1-\left| A_u\left( u\right) \right| -\left| A_u\left(
y_2\right) \right| =\gamma _u-2-\left| A_u\left( y_2\right) \right| 
$$
and therefore%
$$
\begin{array}{ll}
\left| O_u(x,y)\right| -1\geq \mid x 
\overrightarrow{H}y_1\Lambda _u\left( y_1,y_2\right) y_2\overleftarrow{H}%
y\mid -1 \\  \geq \beta +\left| A_u\left( y_1\right) \right| +\left| A_u\left(
y_2\right) \right|  \geq \gamma _u-2+\left| A_u\left( y_1\right) \right|
\geq \gamma _u-1. 
\end{array}
$$

{\bf Case 2.}\quad $u\in \left\{ x,y\right\} .$

Assume w.l.o.g. that $u=x$. Let $x_1\overrightarrow{H}y_1$ be the longest
segment in $x^{+}\overrightarrow{H}y$ with $x_1,y_1\in \Lambda _x$ and $x_2%
\overrightarrow{H}y_2$ be the longest segmnet in $y^{+}\overrightarrow{H}x$
with $y_2\in \Lambda _x$. Putting $\beta =\mid x_1\overrightarrow{H}y\mid +$ 
$\mid x_2\overrightarrow{H}x\mid -2$ we see (by lemma 6) that%
$$
\beta \geq \gamma _x-1-\left| A_x\left( x\right) \right| -\left| A_x\left(
x_1\right) \right| =\gamma _x-2-\left| A_x\left( x_1\right) \right| , 
$$
and therefore%
$$
\begin{array}{ll}
\left| O_x(x,y)\right| -1\geq \mid x 
\overleftarrow{H}x_2\Lambda _x\left( x_2,x_1\right) x_1\overrightarrow{H}%
y\mid -1  \\ \geq \beta +\left| A_x\left( x_1\right) \right| +\left| A_x\left(
x_2\right) \right|  \geq \gamma _x+\left| A_x\left( x_2\right) \right|
-2\geq \gamma _x-1. 
\end{array}
$$

Also, by $(e1)$ and $(e2)$, $\mid x_1\overrightarrow{H}x\mid -1\geq \gamma
_x-1,\mid x_2\overrightarrow{H}y\mid -1\geq \gamma _x-1$ and hence%
$$
\begin{array}{l}
\mid O( 
\stackrel{o}{x},x)\mid -1\geq \mid \stackrel{o}{x}\rho _x\left( x_1\right)
T\left( x_1\right) x_1\overrightarrow{H}x\mid -1\geq \mid x_1\overrightarrow{%
H}x\mid \geq \gamma _x, \\ \mid O(\stackrel{o}{x},y)\mid -1\geq \mid 
\stackrel{o}{x}\rho _x\left( x_2\right) T\left( x_2\right) x_2%
\overrightarrow{H}y\mid -1\geq \mid x_2\overrightarrow{H}y\mid \geq \gamma
_x.\quad \Delta 
\end{array}
$$

${\bf (g4)\quad }$We can suppose $u\not \in \left\{ x,y\right\} $, since
otherwise the arguments are the same. Assume w.l.o.g. that $u\in V(x^{+}%
\overrightarrow{H}y^{-})$. Clearly $\left| O_u(x,y)\right| =\left|
O(x,y)\right| $. In order to prove $(g4.1){\bf -}(g4.4)$, we recall (by the
hypothesis) that $\Lambda _u\subseteq V(x\overrightarrow{H}y)$ and $(%
\stackrel{0}{u},H)\not \in \Delta $.\\

${\bf (g4.1),(g4.2)\quad }$By $(e2)$, $\gamma _u-1=\left| O(x,y)\right|
-1\geq \mid x\overrightarrow{H}y\mid -1\geq \gamma _u-1$ which implies $%
\left| O(x,y)\right| -1=\mid x\overrightarrow{H}y\mid -1=\gamma _u-1$. Using 
$(e3)$, it is easy to see that $(\stackrel{o}{u},x\overrightarrow{H}y)\in
\Delta $, $\gamma _u-1=2(\varphi _u-1)$ and $B_u=\Lambda _u\backslash \{u\}\subseteq U_0$.\\

${\bf (g4.3)\quad }${\bf Case 1.} $z\in U_1.$

{\bf Case 1.1.}\quad $\Lambda _z\not \subseteq V(x\overrightarrow{H}y).$

Choose $w\in \Lambda _z\bigcap V(y^{+}\overrightarrow{H}x^{-}).$ If $%
z=x^{+} $, then%
$$
\left| O(x,y)\right| -1\geq \mid y\overleftarrow{H}z\Lambda _z\left(
z,w\right) w\overrightarrow{H}x\mid -1\geq \mid x\overrightarrow{H}y\mid , 
$$
a contradiction. Otherwise we reach a contradiction by the following way 
$$
\left| O(x,y)\right| -1\geq \mid y\overleftarrow{H}z^{+}\stackrel{0}{u}x^{++}%
\overrightarrow{H}z\Lambda _z\left( z,w\right) w\overrightarrow{H}x\mid
-1\geq \mid x\overrightarrow{H}y\mid . 
$$

{\bf Case 1.2.}\quad $\Lambda _z\subseteq V(x\overrightarrow{H}y).$

Choose $w\in \Lambda _z\backslash \{z\}$. Assume w.l.o.g. that $w\in V(x\overrightarrow{H}%
z^{-})$. Since $z\in U_1$, we have $\left| \Lambda _z\left( z,w\right)
\right| -1\geq 2$ and hence%
$$
\left| O(x,y)\right| -1\geq \mid y\overleftarrow{H}z\Lambda _z\left(
z,w\right) w\overrightarrow{H}z^{-}\stackrel{o}{u}w^{^{\prime }}%
\overleftarrow{H}x\mid -1\geq \mid x\overrightarrow{H}y\mid 
$$
for some $w^{^{\prime }}\in \left\{ w^{-},w^{--}\right\} $, a contradiction.

{\bf Case 2.}\quad $z\in U_0\bigcup U_2.$

If $z\in U_2$, then apply the arguments used in the proof of $(d1)$ (see case
2 and 3). Let $z\in U_0$. If there exists a vertex $w\in \left( \Lambda
_z\backslash \{z\}\right) \backslash\Lambda _u$, then we can reach a contradiction as in case 1.
Otherwise, $\Lambda _z\subseteq \Lambda _u\bigcup \left\{ z\right\} $ and $%
\gamma _z\leq \varphi _u=(\gamma _u+1)/2.$\\

${\bf (g4.4)\quad }$Suppose, to the contrary, that $\Lambda _z\not \subseteq 
$ $V(x\overrightarrow{H}y)$. If $\mid y\overrightarrow{H}x\mid -1=2$, then
clearly $(\stackrel{o}{u},H)\in \Delta $, a contradiction. Let $\mid y%
\overrightarrow{H}x\mid -1\geq 3$. Choose $w\in \Lambda _z\bigcap V(y^{+}%
\overrightarrow{H}x^{-})$. Assume w.l.o.g. that $\mid w\overrightarrow{H}%
x\mid -1\geq 2$. If $z\not \in \Lambda _u$, then by $(g4.3)$, we are done.
Otherwise,%
$$
\left| O(x,y)\right| -1\geq \mid y\overleftarrow{H}z^{++}\stackrel{o}{u}%
x^{++}\overrightarrow{H}z\Lambda _z\left( z,w\right) w\overrightarrow{H}%
x\mid -1\geq \mid x\overrightarrow{H}y\mid , 
$$
a contradiction. So, $(g4.1)$ ${\bf -}$ $(g4.4)$ are proved. A similar proof holds
for $(g4.5)$ ${\bf -}$ $(g4.7)$ when $\Lambda _u\not \subseteq $ $V(x%
\overrightarrow{H}y)$ and $\Lambda _u\not \subseteq $ $V(y\overrightarrow{H}%
x)$. So, the proof of $(g4)$ is completed.\quad$\Delta $\\

${\bf (g6),(g7)\quad }$Let $u=x^{+}$. Choose $v\in \Lambda _u$ so as to
maximize $\mid v\overrightarrow{H}y\mid $. Clearly, $v\in V(y\overrightarrow{%
H}u)$.

{\bf Case 1.}\quad $v=u.$

{\bf Case 1.1}\quad $\left| T\left( u\right) \right| -1\geq 1.$

By $(e1)$ and $(e2)$, $\left| O(x,y)\right| -1\geq \mid u\overrightarrow{H}%
y\mid \geq \gamma _u$. If $\mid O(x,y)\mid -1=\gamma _u$, then by $(e1){\bf ,}%
\left| T\left( u\right) \right| -1\leq 1$ and $\mid u\overrightarrow{H}y\mid
=\gamma _u-1$ which by $(e3)$ holds $(\hat u,u\overrightarrow{H}y)\in \Delta
.$

{\bf Case 1.2}\quad $\left| T\left( u\right) \right| -1=0.$

Clearly, $\left| O(x,y)\right| -1\geq \mid x\overrightarrow{H}y\mid -1\geq
\gamma _u$. If $\mid O(x,y)\mid -1=\gamma _u$, then $\mid x\overrightarrow{H}%
y\mid -1=\gamma _u$ implying that $uw\in E$ for each $w\in V(x%
\overrightarrow{H}y)\backslash \{u\}$, i.e. $(u,x\overrightarrow{H}y)\in \Delta $.

{\bf Case 2.}\quad $v\not =u.$

{\bf Case 2.1}\quad $\left| T\left( u\right) \right| -1\geq 1.$

By $(e1)$ and $(e2)$, $\mid v\overrightarrow{H}y\mid -1\geq \gamma _u-1$ and
hence%
$$
\left| O(x,y)\right| -1\geq \mid y\overleftarrow{H}u\Lambda _u\left(
u,v\right) v\overrightarrow{H}x\mid -1\geq \mid v\overrightarrow{H}y\mid
-1+\left| T\left( u\right) \right| -1\geq \gamma _u+\left| T\left( u\right)
\right| -2. 
$$

If $\left| O(x,y)\right| -1=\gamma _u$, then $\left| T\left( u\right) \right|
-1=1$, $\mid v\overrightarrow{H}y\mid -1=\gamma _u-1$ and by $(e3)$, ($\hat
u,v\overrightarrow{H}y)\in \Delta $.

{\bf Case 2.2}\quad $\left| T\left( u\right) \right| -1=0.$

Clearly, $\left| O(x,y)\right| -1\geq \mid y\overleftarrow{H}u\Lambda
_u\left( u,v\right) v\overrightarrow{H}x\mid -1\geq \mid v\overrightarrow{H}%
y\mid -1\geq \gamma _u.$ If $\left| O(x,y)\right| -1=\gamma _u$, then $\mid v%
\overrightarrow{H}y\mid -1=\gamma _u$ implying that $uw\in E$ for each $w\in
V(v\overrightarrow{H}y)\backslash \{u\}$, i.e. $(u,v\overrightarrow{H}y)\in \Delta $. \quad $\Delta 
$\\

${\bf (g8)\quad }$By $(g4)$, $(\stackrel{o}{x},H)\in \Delta $. Since $%
\left\{ w,y\right\} \subseteq \Lambda _x$, we have $h\geq 6$. If $\mid
O_x(x,y)\mid -1\leq \gamma _z$, then by $(g4.6)$ and $(g4.7){\bf ,}h-2{\bf =}%
\left| O_x(x,y)\right| -1\leq \gamma _z\leq h/2$ implying that $h\leq 4$, a
contradiction. So, $\left| O_x(x,y)\right| -1\geq \gamma _z+1$. By symmetry, 
$\left| O_x(x,w)\right| -1\geq \gamma _z+1$ and the result follows.  \qquad $\Delta $\\

{\bf Proof of lemma 8.} Immediate from definition 3.11.   \qquad $\Delta $\\

{\bf Proof of lemma 9.} By $(d3{\bf )},h\geq \gamma _i+1$ and $h\geq \gamma
_{i+1}+1$ for each $i\in \{1,...,h\}$. In other words, 
\begin{equation}
\label{8}h-1\geq (\gamma _i+\gamma _{i+1})/2\quad \left( i=1,...,h\right) . 
\end{equation}

${\bf (i1)\quad }$By lemma 8, it sufficies to prove $\left| O(x,y)\right|
-1\geq \beta _i$. Assume w.l.o.g. that $i=1$ and $u_{1,}u_2\in V(x^{+}%
\overrightarrow{H}y^{-})$.

{\bf Case 1.}\quad $u_{1,}u_2\in U_0.$

Putting $\Gamma _i=\Phi _i\bigcap V\left( H\right) ~\left( i=1,2\right) $
we see that $\left| \Gamma _i\right| =\varphi _i-b_i^{*}=\gamma _i~$ $\left(
i=1,2\right) $.

{\bf Case 1.1.}\quad $\Gamma _1\bigcup \Gamma _2\subseteq V(x%
\overrightarrow{H}y).$

Clearly $\left| O(x,y)\right| -1\geq \mid x\overrightarrow{H}y\mid -1\geq
\max \left( \left| \Gamma _1\right| ,\left| \Gamma _2\right| \right) \geq
\left( \gamma _1+\gamma _2\right) /2.$

{\bf Case 1.2.}\quad $\Gamma _1\bigcup \Gamma _2\not \subseteq V(x%
\overrightarrow{H}y).$

Assume w.l.o.g. that $\Gamma _1\bigcap (V(y^{+}\overrightarrow{H}x^{-}))\not
=\emptyset $. Let $z_1\overrightarrow{H}z_2$ be the longest segment in $y^{+}%
\overrightarrow{H}x^{-}$ with $z_1,z_2\in \Gamma _1.$

{\bf Case 1.2.1.}\quad $\Gamma _2\bigcap V(y^{+}\overrightarrow{H}%
x^{-})=\emptyset .$

Choose $w\in V(x^{+}\overrightarrow{H}u_1)$ such that $u_2w\in E$ and $\mid x%
\overrightarrow{H}w\mid $ is minimum. Then%
$$
\begin{array}{c}
\left| O(x,y)\right| -1\geq \mid x 
\overrightarrow{H}y\mid -1\geq \mid w\overrightarrow{H}y\mid +\mid x%
\overrightarrow{H}w\mid -2\geq \gamma _2+\mid x\overrightarrow{H}w\mid -2,
\\ \left| O(x,y)\right| -1\geq \mid x\overleftarrow{H}z_1u_1\overleftarrow{H}%
wu_2\overrightarrow{H}y\mid -1\geq \gamma _1-\mid x\overrightarrow{H}w\mid
+2. 
\end{array}
$$

Combining these two inequalities yields $\left| O(x,y)\right| -1\geq (\gamma
_1+\gamma _2)/2.$

{\bf Case 1.2.2.}\quad $\Gamma _2\bigcap V(y^{+}\overrightarrow{H}x^{-})\not
=\emptyset .$

Let $w_1\overrightarrow{H}w_2$ be the longest segment in $y^{+}%
\overrightarrow{H}x^{-}$ with $w_1,w_2\in \Gamma _2.$

{\bf Case 1.2.2.1.}\quad $z_1,w_2\in V(w_1\overrightarrow{H}z_2).$

It follows that $\left| O(x,y)\right| -1\geq \mid x\overrightarrow{H}u_1z_2%
\overleftarrow{H}w_1u_2\overrightarrow{H}y\mid -1\geq $ $\max \left( \gamma
_1,\gamma _2\right) \geq \left( \gamma _1+\gamma _2\right) /2.$

{\bf Case 1.2.2.2.}\quad $z_2,w_1\in V(z_1\overrightarrow{H}w_2).$

Clearly 
$$\left| O(x,y)\right| -1\geq \mid x\overrightarrow{H}u_1z_1%
\overrightarrow{H}w_2u_2\overrightarrow{H}y\mid -1\geq  \max \left( \gamma
_1,\gamma _2\right) \geq \left( \gamma _1+\gamma _2\right) /2.
$$

{\bf Case 1.2.2.3.}\quad Either $w_1,w_2\in V(z_1\overrightarrow{H}z_2)$ or $%
z_1,z_2\in V(w_1\overrightarrow{H}w_2)$.

Assume w.l.o.g. that $w_1,w_2\in V(z_1\overrightarrow{H}z_2)$. If $w_1=z_1$
(resp. $w_2=z_2$), then we can argue as in case 1.2.2.1. (resp. 1.2.2.2.).
Otherwise ($w_1\not =z_1$ and $w_2\not =z_2$),%
$$
\begin{array}{c}
\left| O(x,y)\right| -1\geq \mid x 
\overrightarrow{H}u_1z_1\overrightarrow{H}w_2u_2\overrightarrow{H}y\mid
-1\geq \gamma _2+\mid z_1\overrightarrow{H}w_1\mid -1, \\ \left|
O(x,y)\right| -1\geq \mid x\overrightarrow{H}u_1z_2\overleftarrow{H}w_1u_2%
\overrightarrow{H}y\mid -1\geq \gamma _1-\mid z_1\overrightarrow{H}w_1\mid
+1. 
\end{array}
$$

Combining these two inequalities yields $\left| O(x,y)\right| -1\geq \left(
\gamma _1+\gamma _2\right) /2.$

{\bf Case 2.}\quad $u_1,u_2\in \overline{U}_0.$

By $(g2)$ and $(g3)$, $\left| O(x,y)\right| -1\geq \gamma _i-1~(i=1,2)$. If
either $u_1\in U_{*}$ or $u_2\in U_{*}$, then by $(g1)$ we are done. Let $%
u_1,u_2\in U_1\bigcup U_2.$

{\bf Case 2.1.}\quad Either $\left| O(x,y)\right| -1=\gamma _1-1$ or $\left|
O(x,y)\right| -1=\gamma _2-1$.

Assume w.l.o.g. that $\left| O(x,y)\right| -1=\gamma _1-1.$ Using $(g2)$ we
see that $\left| T\left( u\right) \right| -1=1$ and by $(g4.1)$ and $(g4.3)$%
, $u_2\in U_0\bigcup U_{*},$ a contradiction.

{\bf Case 2.2.}\quad $\left| O(x,y)\right| -1\geq \gamma _1$ and $\left|
O(x,y)\right| -1\geq \gamma _2.$

Clearly, $\left| O(x,y)\right| -1\geq \max \left( \gamma _1,\gamma _2\right)
\geq \left( \gamma _1+\gamma _2\right) /2.$

{\bf Case 3.}\quad $u_1\in U_0,~u_2\in \overline{U}_0.$

By $(g2)$ and $(g3)$, 
\begin{equation}
\label{9}\left| O(x,y)\right| -1\geq \gamma _2-1. 
\end{equation}

{\bf Case 3.1.}\quad $\Phi _1\backslash B_1^{*}\subseteq V(x\overrightarrow{H}y).$

Clearly $\left| O(x,y)\right| -1\geq \mid x\overrightarrow{H}y\mid -1\geq
\gamma _1.$ The result is immediate if either $\left| O(x,y)\right|
-1>\gamma _1$ or $\left| O(x,y)\right| -1>\gamma _2-1$. Thus we can assume $%
\left| O(x,y)\right| -1=\mid x\overrightarrow{H}y\mid -1=\gamma _1=\gamma
_2-1.$ Since $\mid x\overrightarrow{H}y\mid -1=\gamma _1$, we have $\Lambda
_1=V(x\overrightarrow{H}y)\backslash \{u_1\}$. On the other hand, by $(g4.3)$, $\Lambda
_1\subseteq \Lambda _2\bigcup \left\{ u_1\right\} $ , a contradiction.

{\bf Case 3.2.}\quad $\Phi _1\backslash B_1^{*}\not \subseteq V(x\overrightarrow{H}y).$

Let $z_1\overrightarrow{H}z_2$ be the maximal segment in $y^{+}%
\overrightarrow{H}x^{-}$ with $z_1,z_2\in N\left( u_1\right) $.

{\bf Case 3.2.1.}\quad $\Lambda _2\subseteq V(x\overrightarrow{H}y).$

For each $v\in \Lambda _2\bigcap $ $V(x\overrightarrow{H}u_2)$, put $P_v=x%
\overleftarrow{H}z_1u_1\overrightarrow{H}y$ if $v=u_2$ and 
$$
P_v=x\overleftarrow{H}z_1u_1\overleftarrow{H}v\Lambda _2\left( v,u_2\right)
u_2\overrightarrow{H}y 
$$
if $v\not =u_2$. Choose $w_1\in \Lambda _2\bigcap $ $V(x\overrightarrow{H}%
u_2)$ so as to maximize $\mid w_1\overrightarrow{H}u_2\mid $. By $(g2)$ and $%
(g3)$, $\mid w_1\overrightarrow{H}y\mid -1\geq \gamma _2-1.$ If $w_1\not =x$,
then clearly $\mid z_1\overrightarrow{H}y\mid -1\geq \gamma _1$ and hence%
$$
\begin{array}{l}
\left| O(x,y)\right| -1\geq \left| P_{w_1}\right| -1\geq \mid z_1 
\overrightarrow{H}y\mid -\mid x\overrightarrow{H}w_1\mid +2\geq \gamma
_1-\mid x\overrightarrow{H}w_1\mid +3, \\ \left| O(x,y)\right| -1\geq \mid x%
\overrightarrow{H}y\mid -1\geq \mid w_1\overrightarrow{H}y\mid +\mid x%
\overrightarrow{H}w_1\mid -2\geq \gamma _2+\mid x\overrightarrow{H}w_1\mid
-2. 
\end{array}
$$

Combining these two inequalities yields the results. Now let $w_1=x$. Choose 
$w_2\in \Lambda _2\bigcap $ $V(x^{+}\overrightarrow{H}u_2)$ so as to
maximize $\mid w_2\overrightarrow{H}u_2\mid $. Since $H$ is extreme,%
$$
h\geq \mid u_2\Lambda _2\left( u_2,x\right) x\overrightarrow{H}u_1z_2%
\overleftarrow{H}u_2\mid -1 
$$
which implies that 
\begin{equation}
\label{10}
\begin{array}{lll}
\overline{\rho }_2\left( x\right) =\hat u_2 & \Longrightarrow & \mid z_2%
\overrightarrow{H}x\mid -1\geq \left| A_2\left( x\right) \right| +\left|
A_2\left( u_2\right) \right| , 
\end{array}
\end{equation}

\begin{equation}
\label{11}
\begin{array}[t]{lll}
\overline{\rho }_2\left( x\right) =\stackrel{0}{u}_2 & \Longrightarrow & 
\mid z_2\overrightarrow{H}x\mid -1\geq \left| A_2\left( x\right) \right|
+1\geq 2. 
\end{array}
\end{equation}
Observe also, that 
\begin{equation}
\label{12}\left| O(x,y)\right| -1\geq \mid x\overrightarrow{H}y\mid -1\geq
\gamma _1-\mid z_1\overrightarrow{H}z_2\mid . 
\end{equation}

{\bf Claim 9.1.}\quad $\left| O(x,y)\right| -1\geq \gamma _2+\mid z_1%
\overrightarrow{H}z_2\mid .$

{\bf Proof of} {\bf claim 9.1.}\quad Clearly, 
\begin{equation}
\label{13}\left| O(x,y)\right| -1\geq \left| P_{w_2}\right| -1\geq \mid w_2%
\overrightarrow{H}y\mid +\mid z_2\overrightarrow{H}x\mid +\mid z_1%
\overrightarrow{H}z_2\mid -1. 
\end{equation}

By $(f1)$, $\mid w_2\overrightarrow{H}y\mid -1\geq \gamma _2-2-\left|
A_2\left( x\right) \right| $ if $\overline{\rho }_2\left( x\right) =%
\stackrel{o}{u}_2$ and $\mid w_2\overrightarrow{H}y\mid -1\geq \gamma
_2-1-\left| A_2\left( u_2\right) \right| -\left| A_2\left( x\right) \right| $
if $\overline{\rho }_2\left( x\right) =\hat u_2$. Using also (\ref{10}) and (%
\ref{11}), we obtain $\mid w_2\overrightarrow{H}y\mid +\mid z_2%
\overrightarrow{H}x\mid -2\geq \gamma _2-1$, which by (\ref{13}) implies $%
\left| O(x,y)\right| -1\geq \gamma _2+\mid z_1\overrightarrow{H}z_2\mid $. \quad$%
\Delta $

Claim 9.1 with together (\ref{12}) implies the result.

{\bf Case 3.2.2.}\quad $\Lambda _2\not \subseteq V(x\overrightarrow{H}y).$

Let $y_1\overrightarrow{H}y_2$ be the maximal segment in $y^{+}%
\overrightarrow{H}z^{-}$ with $y_1,y_2\in \Lambda _2$.

{\bf Case 3.2.2.1.}\quad Either $z_1,y_2\in V(y_1\overrightarrow{H}z_2)$ or $%
z_2,y_1\in V(z_1\overrightarrow{H}y_2)$.

Assume w.l.o.g that $z_1,y_2\in V(y_1\overrightarrow{H}z_2)$. Then 
$$
\left| O(x,y)\right| -1\geq \mid y\overleftarrow{H}u_2\Lambda _2\left(
u_2,y_1\right) y_1\overrightarrow{H}z_2u_1\overleftarrow{H}x\mid -1\geq
\gamma _1+1, 
$$
and the result follows by (\ref{9}).

{\bf Case 3.2.2.2.}\quad $z_1,z_2\in V(y_1\overrightarrow{H}y_2).$

Apply the arguments in case {\bf \ }3.2.2.1.

{\bf Case 3.2.2.3.}\quad $y_1,y_2\in V(z_1\overrightarrow{H}z_2).$

Putting $\beta =\mid x\overrightarrow{H}y\mid +\mid y_1\overrightarrow{H}%
y_2\mid -2$ and 
$$
\begin{array}{c}
P_1=y 
\overleftarrow{H}u_2\Lambda _2\left( u_2,y_2\right) y_2\overleftarrow{H}%
z_1u_1\overleftarrow{H}x, \\ P_2=y\overleftarrow{H}u_2\Lambda _2\left(
u_2,y_1\right) y_1\overrightarrow{H}z_2u_1\overleftarrow{H}x, 
\end{array}
$$
we obtain 
\begin{equation}
\label{14}\left| O(x,y)\right| -1\geq \left| P_2\right| -1\geq \mid x%
\overrightarrow{H}y\mid +\mid y_1\overrightarrow{H}z_2\mid \geq \gamma
_1-\mid z_1\overrightarrow{H}y_1\mid +2. 
\end{equation}

{\bf Claim 9.2.}\quad $\left| O(x,y)\right| -1\geq \gamma _2+\mid z_1%
\overrightarrow{H}y_1\mid -1.$

{\bf Proof of} {\bf claim 9.2.}\quad If $\overline{\rho }_2\left( y_2\right)
=\hat u_2$, then by $(f1)$, $\beta \geq \gamma _2-1-\left| A_2\left(
u_2\right) \right| -\left| A_2\left( y_2\right) \right| $ and%
$$
\begin{array}{c}
\left| O(x,y)\right| -1\geq \left| P_1\right| -1\geq \beta +\left| A_2\left(
u_2\right) \right| +\mid A_2\left( y_2\right) \mid +\mid z_1 
\overrightarrow{H}y_1\mid -1 \\ \geq \gamma _2+\mid z_1\overrightarrow{H}%
y_1\mid -2. 
\end{array}
$$

Otherwise $(\overline{\rho }_2(y_2)=\stackrel{o}{u}_2)$, $\beta \geq \gamma
_2-2-\left| A_2\left( y_2\right) \right| =\gamma _2-3$ and 
$$
\left| O(x,y)\right| -1\geq \left| P_1\right| -1\geq \beta +\mid z_1%
\overrightarrow{H}y_1\mid +1\geq \gamma _2+\mid z_1\overrightarrow{H}y_1\mid
-2.  \quad \Delta 
$$

Claim 9.2 together with (\ref{14}) implies the result.\quad$\Delta $\\

${\bf (i2)\quad }$By $(d1)$,~$h\geq 2\gamma _x$ imlying that $\left|
O(x,y)\right| -1\geq h/2\geq \gamma _x$. Also, $\left| O(x,y)\right| -1\geq
\gamma _z$ by $(g6)$. Using lemma 8, $\left| \Omega \left( x,y\right)
\right| -1\geq \left| O(x,y)\right| -1\geq \left( \gamma _x+\gamma _z\right)
/2.  \quad \Delta $\\

${\bf (i3)\quad }$By $(g1)$, $\left| O(x,y)\right| -1\geq \gamma _x+1$ and
by $(g6)$, $\left| O(x,y)\right| -1\geq \gamma _z$. Using lemma 8, we obtain
the result immediately.  \quad $\Delta $\\

${\bf (i4)\quad }${\bf Claim 9.3.\quad }$\max \left( \left| O_x(x,y)\right|
-1,\left| O_x(x,w)\right| -1\right) \geq \left( \gamma _x+\gamma _z\right)
/2.$

{\bf Proof of} {\bf claim 9.3.}\quad By $(g6){\bf ,~}\min (\left|
O(x,y)\right| -1,\left| O(x,w)\right| -1)\geq \gamma _z$. If either $\left|
O_x(x,y)\right| -1\geq \gamma _x$ or $\left| O_x(x,w)\right| -1\geq \gamma
_x $, then clearly we are done. Otherwise, by $(g3){\bf ,}\left|
O_x(x,y)\right| -1=\left| O_x(x,w)\right| -1=\gamma _x-1$ and the result
holds by $(g8)$ and lemma 8.  \quad $\Delta $\\

{\bf Claim 9.4.}\quad $\min (\mid O(\stackrel{o}{x},y)\mid -1,\mid O(%
\stackrel{o}{x},w)\mid -1)\geq \left( \gamma _x+\gamma _z+1\right) /2.$

{\bf Proof of} {\bf claim 9.4.}\quad By $(g5)$ and $(g6){\bf ,}\mid O(%
\stackrel{o}{x},y)\mid -1\geq \gamma _x{\bf ~}$ and $\left| O(x,y)\right|
-1\geq \gamma _z$, respectively. Since $V(O(x,y))\bigcap \{\stackrel{o}{x}%
\}=\emptyset $ (by definition 2.10), $\left| O(\stackrel{o}{x},y)\right|
-1\geq \left| O(x,y)\right| \geq \gamma _z+1$, implying that $\mid O(%
\stackrel{o}{x},y)\mid -1\geq (\gamma _x+\gamma _z+1)/2$. Analogously, $\mid
O(\stackrel{o}{x},w)\mid -1\geq (\gamma _x+\gamma _z+1)/2$.  \quad  $\Delta $\\

{\bf Claim 9.5.}\quad $\mid O(\stackrel{o}{x},x)\mid -1\geq (\gamma
_x+\gamma _z+1)/2$.

{\bf Proof of} {\bf claim 9.5.}\quad Let $v\in \Lambda _x\backslash \{x\}$. By $(g5)$ and $%
(g6)$, $\mid O(\stackrel{o}{x},x)\mid -1\geq \gamma _x$ and $\mid O(x,v)\mid
-1\geq \gamma _z$, respectively. Hence%
$$
\mid O(\stackrel{o}{x},x)\mid -1\geq \left| O(x,v)\right| +\mid vT\left(
v\right) \stackrel{\wedge }{v}\stackrel{o}{x}\mid -2\geq \gamma _z+1 
$$
which implies $\mid O(\stackrel{o}{x},x)\mid -1\geq \left( \gamma _x+\gamma
_z+1\right) /2.   \quad\Delta $

The result holds from claims 9.3-9.5 and lemma 8.$ \quad \Delta $\\

${\bf (i5)\quad }$By $(g6)$, $\left| O(x,y)\right| -1\geq \gamma _z$ and $%
\left| O(z,w)\right| -1\geq \gamma _x$ and the result follows from lemma 8.  \quad $%
\Delta $\\

${\bf (i6)\quad }$By $(g6)$, $\left| O(x,y)\right| -1\geq \max \left( \gamma
_{x^{+}},\gamma _{x^{-}}\right) .$ If $\left| T\left( x\right) \right|
-1\geq 2$, then by $(g2){\bf ,}\left| O(x,y)\right| -1\geq \gamma _x$ and the
result holds immediately. Thus we can assume $\left| T\left( x\right)
\right| -1=1.$ Put $z=x^{+}$ and $w=x^{-}$. By $(g3)$ and $(g6){\bf ,~}%
\left| O_x(x,y)\right| -1\geq \max \left( \gamma _x-1,\gamma _z,\gamma
_w\right) .$ If ${\bf ~}\left| O_x(x,y)\right| -1\geq \min \left( \gamma
_x,\gamma _z+1,\gamma _w+1\right) $, then clearly we are done. Now let $%
\left| O_x(x,y)\right| -1=\gamma _x-1=\gamma _z=\gamma _w$. Since $u\not \in
U_{*}$ (by $(g1)$), we have by $(g4.3)$ and $(g4.7){\bf ,}\gamma _z\leq
\left( \gamma _x+1\right) /2=\left( \gamma _z+2\right) /2$ implying that $%
\gamma _z\leq 2$ and $\left| O_x(x,y)\right| -1\leq 2$. It means that $h\leq
4$. Recalling also $(g4.1)$ and $(g4.5)$, we conclude that $h=4$, a
contradiction. \quad  $\Delta $\\

${\bf (i7)\quad }$If either $\left| O_x(x,y)\right| -1=\gamma _x-1$ or $%
\left| O_y(x,y)\right| -1=\gamma _y-1$, say $\left| O_x(x,y)\right|
-1=\gamma _x-1$, then by $(g1){\bf -}(g3){\bf ,~}\left| T\left( x\right)
\right| -1=1$. By $(g4)$, either $(\stackrel{o}{x},x\overrightarrow{H}y)\in
\Delta $ or $(\stackrel{o}{x},y\overrightarrow{H}x)\in \Delta $ or $(%
\stackrel{o}{x},H)\in \Delta $. This implies by $(g4.2)$ and $(g4.6)$ that $%
y\in U_0$, a contradiction. Thus $\left| O_x(x,y)\right| -1\geq \gamma _x$
and $\left| O_y(x,y)\right| -1\geq \gamma _y$. Using $(g6)$ with lemma 8, we
obtain $\left| \Omega (x,y)\right| -1\geq \left( \gamma _x+\gamma _z\right)
/2$ for each $z\in \left\{ x^{+},x^{-}\right\} $ and $\left| \Omega
(x,y)\right| -1\geq \left( \gamma _y+\gamma _w\right) /2$ for each $w\in
\left\{ y^{+},y^{-}\right\} $. Then the result follows from $(i1)$.\quad$\Delta $\\

${\bf (i8)\quad }$Observing that $\left| O(x,y)\right| -1\geq \mid y%
\overrightarrow{H}x\mid -1=h-1$, we obtain the result from (\ref{8})
immediately.\quad$\Delta $\\

${\bf (i9)\quad }$Put $z=x^{+}$. We can assume $h\geq 4$, since otherwise
the result holds from $(i8)$. By $(d3){\bf ,}\left| O(x,y)\right| -1\geq
h-2\geq \gamma _x-1$. If $\left| O(x,y)\right| -1\geq \gamma _z+1$, then
clearly we are done. Let $\left| O(x,y)\right| -1=\gamma _z$. By $(g7){\bf ,~%
}\left| T\left( z\right) \right| -1\leq 1$. If $\Lambda _z\bigcap V(y^{+}%
\overrightarrow{H}x^{-})\not =\emptyset $, then by ${\bf (g7),}~\stackrel{%
\wedge }{z}x^{-}\in E$. Hence $\left| O(x,y)\right| -1\geq h-1$ and by (\ref
{8}) we are done. Now let $\Lambda _z\subseteq V(x\overrightarrow{H}y)$. It
means that $\gamma _z=2$ and $\left| O(x,y)\right| -1=2$. But then $h=4$, a
contradiction.      \qquad  $\Delta $\\

\noindent{\bf Proof of the theorem.}\quad Let $G$ be a graph, $C$ be a longest cycle in $G$
and $H=u_1...u_hu_1$ a longest cycle in $G\backslash C$ with a maximal $HC-$extension $%
T$. Putting $U_{*}=\left\{ v_1^{*},...,v_r^{*}\right\} $ and using
definition 3.3, we let for each $i\in \{1,...,r\}$,%
$$
\begin{array}{c}
\Theta ( 
\overleftarrow{T}(v_i^{*}),V_{neut},V_{fin}^{\left( i\right) })=(P_0^{\left(
i\right) },...,P_{\pi \left( i\right) }^{\left( i\right) }), \\ R_i=<(V(%
\stackrel{\wedge }{v}_i^{*}\overleftarrow{T}(v_i^{*})z_{\pi \left( i\right)
}^{\left( i\right) })\bigcup\bigcup^{\pi(i)}_{j=0}
V(P_j^{\left( i\right) }))\backslash \{z_{\pi \left( i\right) }^{\left(
i\right) }\}>, 
\end{array}
$$
where $V_{neut}=V\backslash(V\left( C\right) \bigcup V\left( T\right) )$ and $%
V_{fin}^{\left( i\right) }=V\left( T\right) \backslash V(T(v_i^{*})).$ Since $c\geq
\delta +1\geq \kappa+1$, 
for each $i\in \{1,...,r\}$ there are $\kappa-1$ internally disjoint paths $%
E_i^{\left( 1\right) },...,E_i^{\left( \kappa-1\right) }$ in $\left( \kappa-1\right) -$%
connected graph $G\backslash \{z_{\pi \left( i\right) }^{\left( i\right) }\}$, starting at 
$R_i$, passing through $V_{neut}$ and terminating on $C$ at $\kappa-1$ different
vertices. Let $E_j^{\left( a\right) }$ has a vertex $v$ in common with $%
E_e^{\left( b\right) }$ for some $a,b\in \{1,...,\kappa-1\}$ and $j,e\in 
\{1,...,r\}~(j\not =e)$. If $v\not \in V\left( C\right) $, then there is a
path starting in $R_j$, passing through $V_{neut}$ and terminating in $R_e$,
contradicting the fact that $v_j^{*},v_e^{*}\in U_{*}$. So, $v\in V\left(
C\right) .$ Choose vertex-disjoint paths 
$E_1^{\left( i_1\right)},...,E_t^{\left( i_t\right) }$ 
$\left( i_j\in \{1,...,\kappa-1\}\right)$ 
for each $j\in \{1,...,t\}$ 
so as to maximize $t$ and put $E_j^{\left(
i_j\right) }=x_j\overrightarrow{E}_j^{\left( i_j\right) }w_j^{*}~\left(
j=1,...,t\right) ,$ where $x_j\in V\left( R_j\right) $ and $w_j^{*}\in
V\left( C\right) $. It is easy to see that $t\geq \min \left( r,\kappa-1\right) $%
. By $(a2)$, for each $j\in \{1,...,t\}$ there is an $\left(
x_j,v_j^{*}\right) -$path $F_j^{\left( i_j\right) }$ passing through $%
V\left( R_j\right) \bigcup V(T(v_j^{*}))$ and having length at least $%
\varphi _{v_j^{*}}$. Denoting
$$
E_j^{*}=v_j^{*}F_j^{\left( i_j\right)
}x_jE_j^{\left( i_j\right) }w_j^{*}\quad (~j=1,...,t),
$$
we see that $E_1^{*},...,E_t^{*}$ are vertex disjoint $\left( H,C\right) -$paths with 
$\left| E_i^{*}\right| -1\geq \varphi _{v_i^{*}}+1\quad (i=1,...,t)$.

{\bf Case 1.}\quad $\kappa\geq 4,~h\geq 5.$

{\bf Case 1.1.}\quad $r\geq \kappa.$

It follows that $t\geq \kappa-1$. Let $\xi _1,...,\xi _t$ be the elements of $%
\left\{ w_1^{*},...,w_t^{*}\right\} $ occuring on $C$ in a consequtive order.

{\bf Case 1.1.1.}\quad $t\geq \kappa.$

Assume w.l.o.g. that $\varphi _{v_1^{*}}\geq ...\geq \varphi _{v_r^{*}}$.
Since $r\geq \kappa$, we have%
$$
\frac 1\kappa\sum^\kappa_{i=1}\varphi _{v_i^{*}}\geq \frac 1r%
\sum^r_{i=1}\varphi _{v_i^{*}} 
$$
implying that 
\begin{equation}
\label{15}\sum^\kappa_{i=1}\varphi _{v_i^{*}}\geq \frac
\kappa r\sum^r_{i=1}\varphi _{v_i^{*}}\geq \frac \kappa h%
\sum^r_{i=1}\varphi _{v_i^{*}}. 
\end{equation}

By $(i1)$ and $(i3)$, $\left| \Omega \left( v_a^{*},v_b^{*}\right) \right|
-1\geq \beta _i$ for each $a,b\in \{1,...,t\}$ and $i\in \{1,...,h\}$.
Hence%
$$
\left| \Omega \left( v_a^{*},v_b^{*}\right) \right| -1\geq \frac 1h%
\sum^h_{i=1}\beta _i=\mu \left( T\right) . 
$$

Then for each $i,j\in \{1,...,t\}$, 
\begin{equation}
\label{16}\mid w_i^{*}\overrightarrow{C}w_j^{*}\mid -1\geq \left|
E_i^{*}\right| -1+\left| E_j^{*}\right| -1+\left| \Omega \left(
v_i^{*},v_j^{*}\right) \right| -1\geq \varphi _{v_i^{*}}+\varphi
_{v_j^{*}}+2+\mu \left( T\right) . 
\end{equation}

Using (\ref{15}),(\ref{16}) and recalling that $t\geq \kappa$, we obtain%
$$
\begin{array}{lll}
c= 
\sum^t_{i=1}(\mid \xi _i\overrightarrow{C}\xi
_{i+1}\mid -1)\geq 2\sum^t_{i=1}\varphi
_{v_i^{*}}+2t+t\mu \left( T\right) \\ \geq 
\sum^t_{i=1}\varphi _{v_i^{*}}+2t+t\mu \left(
T\right) \geq \frac \kappa h\sum^r_{i=1}\varphi
_{v_i^{*}}+\kappa\mu \left( T\right) +2\kappa \\ \geq \frac \kappa h
(\sum^r_{i=1}\varphi _{v_i^{*}}+\sum^h_{i=1}
\varphi _{_i}^{^{\prime }}+2h)=\frac \kappa h(\sum^h_{i=1}
\varphi _{_i}^{}+2h), 
\end{array}
$$
where $\xi _{t+1}=\xi _1$. It follows that $\sum_{i=1}^h\varphi _{_i}\leq
h\left( c/\kappa-2\right) .$ Since $\varphi _i+\psi _i=d\left( u_i\right) \geq
\delta \, \left( i=1,...,h\right) $, we have
$$
\sum^h_{i=1}\psi _i\geq h\delta -
\sum^h_{i=1}\varphi _{_i}^{}\geq h\delta -ch/\kappa+2h. 
$$

In particular, $\max_i\psi _i\geq \delta -c/\kappa+2$. Using lemma
3, we obtain%
$$
c\geq \sum^h_{i=1}\psi_i+ \max_i\psi
_i\geq h\delta -ch/\kappa+2h+\delta -c/\kappa+2, 
$$
and the result follows immediately.

{\bf Case 1.1.2.}\quad $t=\kappa-1.$

Observe that $E_\kappa^{\left( i\right) }$ terminates in $\left\{
w_1^{*},...,w_{\kappa-1}^{*}\right\} $ for each $i\in \{1,...,\kappa-1\}$, since
otherwise $E_1^{*},...,E_{\kappa-1}^{*},E_\kappa^{\left( j\right) }$ contradict the
maximality of $t$ for some $j\in \{1,...,\kappa-1\}$. By the same arguments, $%
E_j^{\left( i\right) }$ terminates in $\left\{
w_1^{*},...,w_{\kappa-1}^{*}\right\} $ for each $i\in \{1,...,\kappa-1\}$ and $j\in 
\{1,...,\kappa\}$. Then there is a path $E=vE\xi
_{t+1}$ starting in%
$$
<(V\left( T\right) \bigcup \bigcup^\kappa_{i=1}
R_i\bigcup \bigcup^\kappa_{j=1}
\bigcup^{\kappa-1}_{i=1}V(E_j^{\left( i\right) }))\backslash\left\{
w_1^{*},...,w_{\kappa-1}^{*}\right\} > 
$$
and terminating in $C\backslash\{w_1^{*},...,w_{\kappa-1}^{*}\}$. Assume w.l.o.g. that $%
\xi _1,...,\xi _{t+1}$ occurs on $\overrightarrow{C}$ in a consequtive order.
Then it is easy to see that 
$$
c=\sum^{t+1}_{i=1}(\mid \xi _i\overrightarrow{C}\xi
_{i+1}\mid -1)\geq \sum^\kappa_{i=1}\varphi
_{v_i^{*}}+2\kappa+\kappa\mu \left( T\right) 
$$
where $\xi _{t+2}=\xi _1$. Further, we can argue exactly as in case 1.1.1.

{\bf Case 1.2.}\quad $r\leq \kappa-1.$

It follows that $t=r$. There are $\kappa$
vertex-disjoint $(H,C)-$paths $E_i=v_iE_iw_i\quad \left( i=1,...,\kappa\right) $.
Assume w.l.o.g. that $w_1,...,w_\kappa$ occurs on $\overrightarrow{C}$ in a
consequtive order. Put%
$$
W=\left\{ w_1,...,w_\kappa\right\} ,\quad W^{*}=\left\{
w_1^{*},...,w_r^{*}\right\} . 
$$

Let $a,b\in \{1,...,\kappa\}$. Denote
$$
W^{*}\left( a,b\right) =W^{*}\bigcap V(w_a\overrightarrow{C}w_b) 
$$

We will say that $w_a\overrightarrow{C}w_b$ is a suitable segment if

$$
\mid w_a\overrightarrow{C}w_b\mid -1\geq \sum_{v_i^*\in W^*\left(
a,b\right)}\varphi_{v_i^*}+2\left( b-a\right)
+\sum^{b-a}_{i=1}\left( \left| \Omega \left( \overline{v}_{a+i-1},
\overline{\overline{v}}_{a+i}\right) \right| -1\right), 
$$

where $\overline{v}_j,\overline{\overline{v}}_j\in \left\{ v_j\right\}
\bigcup \overline{U}_0\ (j=1,...,\kappa)$.\\

{\bf Claim 1.}\quad Let $i\in \{1,...,\kappa\}$. If either $\left| W^{*}\left(
i,i+1\right) \right| \not =1$ or $\left| W^{*}\left( i,i+1\right) \right| =1$
and \linebreak
$W^{*}\left( i,i+1\right) \bigcap \left\{ w_i,w_{i+1}\right\}
=\emptyset $, then $w_i\overrightarrow{C}w_{i+1}$ is suitable.

{\bf Proof of claim 1.\quad Case a1}\quad $\left| W^{*}\left( i,i+1\right)
\right| =0.$

Let $T_{tr}\left( E_i,E_{i+1}\right) =(E_i^{^{\prime }},E_{i+1}^{^{\prime
}}) $ and $T_{tr}\left( v_i,v_{i+1}\right) =\left( \overline{v}_i,\overline{%
\overline{v}}_{i+1}\right) $. Then $w_i\overrightarrow{C}w_{i+1}$ is
suitable, since by $(a1){\bf ,}$%
$$
\mid w_i\overrightarrow{C}w_{i+1}\mid -1\geq 2+\left( \left| \Omega \left( 
\overline{v}_i,\overline{\overline{v}}_{i+1}\right) \right| -1\right) . 
$$

{\bf Case a2.}\quad $\left| W^{*}\left( i,i+1\right) \right| \geq 2.$

Let $E,F$ be any two elements of $\left\{ E_1^{*},...,E_r^{*}\right\} $ with 
$E=xEv,F=yFw$ for some $v,w\in W_i^{*}$. Since $T_{tr}\left( E,F\right)
=\left( E,F\right) $ and $\left\{ x,y\right\} \subseteq \overline{U}_0$, we
have by $(a1){\bf ,}\mid v\overrightarrow{C}w\mid -1\geq \varphi _x+\varphi
_y+2+\left( \left| \Omega \left( x,y\right) \right| -1\right) $ implying
that $w_i\overrightarrow{C}w_{i+1}$ is suitable.

{\bf Case a3.}\quad $\left| W^{*}\left( i,i+1\right) \right| =1.$

Assume w.l.o.g that $W^{*}\left( i,i+1\right) =\left\{ w_1^{*}\right\} $. If
either $E_i$ or $E_{i+1}$ (say $E_i$) has no vertex in common with $E_1^{*}$,
then using transformation $T_{tr}\left( E_i,E_1^{*}\right) =(E_i^{^{\prime
}},E_1^{*})$, we obtain by $(a1)$,%
$$
\mid w_i\overrightarrow{C}w_{i+1}\mid -1\geq \mid w_i\overrightarrow{C}%
w_1^{*}\mid -1\geq \varphi _{v_1^{*}}+2+\left( \left| \Omega \left( 
\overline{v}_i,\overline{\overline{v}}_{i+1}\right) \right| -1\right) 
$$
for some appropriate $\overline{v}_i\in \left\{ v_i\right\} \bigcup 
\overline{U}_0$ and $\overline{\overline{v}}_{i+1}=v_1^{*}.$ It means that $%
w_i\overrightarrow{C}w_{i+1}$ is suitable. Now let both $E_i$ and $E_{i+1}$
have vertices in common with $E_1^{*}$. Walking along $E_1^{*}$ from $w_1^{*}$
to $v_1^{*}$ we stop at the first vertex $v\in V\left( E_i\right) \bigcup
V\left( E_{i+1}\right) .$ Assume.w.l.o.g. that $v\in V\left( E_{i+1}\right) $%
. Putting $E_{i+1}^{^{\prime }}=w_1^{*}E_1^{*}vE_{i+1}v_{i+1}$ and $%
T_{tr}(E_i,E_{i+1}^{^{\prime }})=(E_i^{^{\prime }},E_{i+1}^{^{\prime \prime
}})$, we see by $(a1)$ that for some appropriate $\overline{v}_i\in \left\{
v_i\right\} \bigcup \overline{U}_0$ and $\overline{\overline{v}}_{i+1}\in
\left\{ v_{i+1}\right\} \bigcup \overline{U}_0,$%
$$
\mid w_i\overrightarrow{C}w_{i+1}\mid -1\geq \mid w_i\overrightarrow{C}%
w_1^{*}\mid -1\geq 2+\varphi _{v_1^{*}}+\left( \left| \Omega \left( 
\overline{v}_i,\overline{\overline{v}}_{i+1}\right) \right| -1\right) . 
$$

So, again $w_i\overrightarrow{C}w_{i+1}$ is suitable.\quad$\Delta $\\

{\bf Claim 2.}\quad If $w_a\overrightarrow{C}w_b$ and $w_b\overrightarrow{C}%
w_e$ are suitable segments, then $\mid w_a\overrightarrow{C}w_e\mid $ is
suitable as well.

{\bf Proof of claim 2.}\quad Immediate from the definition.\quad$\Delta $\\

{\bf Claim 3.}\quad Let $w_a\overrightarrow{C}w_b$ is a suitable segment. If 
$w_b\overrightarrow{C}w_{b+1}$ is not suitable and \linebreak
$W^{*}\left( b,b+1\right)
=\{w_b\}$, then $w_a\overrightarrow{C}w_{b+1}$ is suitable.

{\bf Proof of claim 3.}\quad Immediate from the definition.\quad$\Delta $\\

{\bf Claim 4.}\quad Let $i\in \{1,...,\kappa\}.$ If $W^{*}\left( i,i+1\right)
\subseteq \left\{ w_i,w_{i+1}\right\} $ and $\left| W^{*}\left( i,i+1\right)
\right|\\ =1$ (say $W^{*}\left( i,i+1\right) =\left\{ w_i\right\} $), then $%
w_{i-1}\overrightarrow{C}w_{i+1}$ is suitable.

{\bf Proof of claim 4.}\quad Assume w.l.o.g. that $W^{*}\left( i,i+1\right)
=\left\{ w_1^{*}\right\} $, i.e. $w_1^{*}=w_i$. If $\left| W^{*}\left(
i-1,i\right) \right| \geq 2$, then by claims 1 and 3, $w_{i-1}\overrightarrow{%
C}w_{i+1}$ is suitable. Let $\left| W^{*}\left( i-1,i\right) \right| =1$,
i.e $W^{*}\left( i-1,i\right) =\left\{ w_1^{*}\right\} $. If either $E_{i-1}$
or $E_{i+1}$ (say $E_{i-1}$) has no vertices in common with $E_1^{*}$, then
using transformations $T_{tr}\left( E_{i-1},E_1^{*}\right) $ and $%
T_{tr}\left( E_i,E_{i+1}\right) $, we see that $w_{i-1}\overrightarrow{C}w_i$
is suitable and by claim 3, $w_{i-1}\overrightarrow{C}w_{i+1}$ is suitable
as well. Now let both $E_{i-1}$ and $E_{i+1}$ have vertices in common with $%
E_1^{*}$. Walking along $E_1^{*}$ from $w_1^{*}$ to $v_1^{*}$ we stop at the
first vertex $v\in V\left( E_{i-1}\right) \bigcup V\left( E_{i+1}\right) .$
Assume w.l.o.g. that $v\in V\left( E_{i+1}\right) $. If $v=\stackrel{o}{v}%
_1^{*}$, i.e. $v_{i-1}=v_1^{*}$, then using $T_{tr}\left(
E_{i-1},w_iE_1^{*}vE_{i+1}v_{i+1}\right) $ and $T_{tr}\left(
E_i,w_{i+1}E_{i+1}vv_1^{*}\right) $ we see that $w_{i-1}\overrightarrow{C}%
w_i $ is suitable. By claim 3, $w_{i-1}\overrightarrow{C}w_{i+1}$ is
suitable as well. Finally, if $v\not =\stackrel{o}{v}_i^{*}$ (i. e. $v_{i+1}%
\not \in U_0)$, then using $T_{tr}\left( E_{i-1},w_iE_1^{*}vE_iv_{i+1}\right) 
$ and $T_{tr}\left( E_i,E_{i+1}\right) $, we see that $w_{i-1}%
\overrightarrow{C}w_i$ is suitable, implying by claim 3 that $w_{i-1}%
\overrightarrow{C}w_{i+1}$ is suitable as well.\quad$\Delta $\\

{\bf Claim 5.}\quad For appropriate $\overline{v}_i,\overline{\overline{v}}%
_i\in \left\{ v_i\right\} \bigcup \overline{U}_0$,%
$$
c=\sum^\kappa_{i=1}(\mid w_i\overrightarrow{C}%
w_{i+1}\mid -1)\geq \sum^r_{i=1}\varphi
_{v_i^{*}}+2\kappa+\sum^\kappa_{i=1}\left( \left| \Omega
\left( \overline{v}_i,\overline{\overline{v}}_{i+1}\right) \right| -1\right)
. 
$$

{\bf Proof of claim 5.}\quad Suppose not. Let $i\in \{1,...,\kappa\}$. If $w_i%
\overrightarrow{C}w_{i+1}$ is not suitable, then by claims 1 and 4, either $%
w_{i-1}\overrightarrow{C}w_{i+1}$ or $w_i\overrightarrow{C}w_{i+2}$ is
suitable. Thus there exist some suitable segment on $C$ and let $w_a%
\overrightarrow{C}w_b$ be the longest one for some $a,b\in \{1,...,\kappa\}%
\quad \left( a\not =b\right) $. If $w_b\overrightarrow{C}w_{b+1}$ is
suitable, then by claim 2, $w_a\overrightarrow{C}w_{b+1}$ is suitable as
well, a contradiction. Otherwise, by claims 3 and 4, $w_b\overrightarrow{C}%
w_{b+2}$ is suitable and hence (by claim 2) $w_a\overrightarrow{C}w_{b+2}$
is suitable as well, a contradiction.\quad$\Delta $\\

{\bf Claim 6.}\ If $\kappa\geq 4$ and $h\geq 5$, then for appropriate $\overline{v}%
_i,\overline{\overline{v}}_i\in \left\{ v_i\right\} \bigcup \overline{U}_0$ 
$\left( i=1,...,\kappa\right) $,%
$$
\sum^\kappa_{i=1}\left( \left| \Omega \left( \overline{v}%
_i,\overline{\overline{v}}_{i+1}\right) \right| -1\right) \geq \kappa\mu \left(
T\right) . 
$$

{\bf Proof of claim 6.}\quad Assume w.l.o.g that $\beta _1=
\max_i\left\{ \beta _i\right\} $. Put%
$$
\begin{array}{ll}
A_0=\left\{ i\left| \left| \Omega \left( 
\overline{v}_i,\overline{\overline{v}}_{i+1}\right) \right| -1\right. \geq
\beta _1\right\} ,\quad A_1=\{1,...,\kappa\}\backslash A_0, \\ A_{1j}=\{i\in A_1\left|
u_j\in \left\{ \overline{v}_i,\overline{\overline{v}}_{i+1}\right\} \right.
\}\quad \left( j=1,2\right) . 
\end{array}
$$

We can assume that $A_1\not =\emptyset $, since otherwise by $(i1)$, 

$$%
\sum_{i=1}^\kappa\left( \left| \Omega \left( \overline{v}_i,\overline{\overline{v}%
}_{i+1}\right) \right| -1\right) \ \geq \kappa\beta _1\geq \kappa\mu \left( T\right). 
$$
If $\left\{ \overline{v}_i,\overline{\overline{v}}_{i+1}\right\} =\left\{
u_1,u_2\right\} $, then by $(i8){\bf ,\ }i\in A_0$. It means that $%
A_{11}\bigcap A_{12}=\emptyset .$ On the other hand, by $(i5)$, either $%
A_{11}=\emptyset $ or $A_{12}=\emptyset $. Assume w.l.o.g. that $%
A_{12}=\emptyset ,$ i.e. $A_1=$ $A_{11}.$

{\bf Case b1.}$\quad \left| A_{11}\right| \geq 4.$

Recalling definition 3.9, it is not hard to see that there are at least two
paths among $E_1,...,E_\kappa$ having vertices in common with $V\left( T\left(
u_1\right) \right) \backslash \{u_1\}$, i.e. $\left| T\left( u_1\right) \right| -1\geq 2$.
By $(i2)$, $A_1=\emptyset $ , a contradiction.

{\bf Case b2.}$\quad \left| A_{11}\right| =3.$

If follows that at least one of the paths $E_1,...,E_\kappa$ has a vertex
in common with $V\left( T\left( u_1\right) \right) \backslash \{u_1\}$, i.e. $\left| T\left(
u_1\right) \right| -1\geq 1$. Clearly $\left| T\left( u_1\right) \right|
-1=1,$ since otherwise $\left| A_{11}\right| \geq 4$. Assume w.l.o.g. that $%
A_1=\left\{ 1,2,3\right\} $ and $\overline{v}_1=$ $\overline{v}_2=$ $%
\overline{v}_3=$ $u_1$. If $\overline{\overline{v}}_2=\overline{\overline{v}}%
_3=\overline{\overline{v}}_4$, then clearly $\overline{v}_1,\overline{%
\overline{v}}_2\in \overline{U}_0$ and by $(i7){\bf ,}\left| \Omega \left( 
\overline{v}_1,\overline{\overline{v}}_2\right) \right| -1\geq \beta _1$, a
contradiction. Otherwise, by $(i4){\bf ,}\left| \Omega \left( \overline{v}_i,%
\overline{\overline{v}}_{i+1}\right) \right| -1\geq \beta _1$ for some $i\in
A_1$, again a contradiction. So, $1\leq \left| A_{11}\right| \leq 2.$

{\bf Case b3.}$\quad \left| A_{11}\right| =2.$

{\bf Case b3.1.}$\quad h\geq 8.$

Let $i\in A_1.$ Assume w.l.o.g. that $\overline{v}_i=u_1,\overline{\overline{%
v}}_{i+1}=u_s$ for some $s\in \{1,...,h\}$. By $(i1)$, $\left| \Omega
\left( \overline{v}_i,\overline{\overline{v}}_{i+1}\right) \right| -1\geq
\beta _j$ for each $j\in \{1,...,h\}\backslash\left\{ 1,h,s-1,s\right\} $. Since $%
h\geq 8$, there are at least four pairwise different integers $%
f_1,f_2,f_3,f_4$ in $\{1,...,h\}\backslash$ $\left\{ 1,h,s-1,s\right\} $. By $(i1)%
{\bf ,}\left| \Omega \left( \overline{v}_i,\overline{\overline{v}}%
_{i+1}\right) \right| -1\geq \beta _{f_j}\quad \left( j=1,2,3,4\right) .$ So,%
$$
\mbox{if} \enskip i\in A_1, \enskip\mbox{then}\enskip \left| \Omega \left( \overline{v}_i,\overline{%
\overline{v}}_{i+1}\right) \right| -1\geq \frac 1h(
\sum^h_{i=1}\beta _i-\beta _1-\beta _h-\beta _{s-1}-\beta _s+%
\sum^4_{i=1}\beta _{f_i}).~ 
$$

On the other hand, 
\begin{equation}
\label{17}\mbox{if}\enskip i\in A_0, \enskip\mbox{then}\enskip  \left| \Omega \left( \overline{v}_i,%
\overline{\overline{v}}_{i+1}\right) \right| -1\geq \beta _1=\frac 1h
(\sum^h_{i=1}\beta _i-
\sum ^h_{i=1}\beta _i+h\beta _1). 
\end{equation}

Since $h\beta _1-\beta _1-\beta _h-\beta _{s-1}-\beta _s\geq
\sum_{i=1}^h\beta _i-\sum_{i=1}^4\beta _{f_i}$, we have%
$$
\mbox{if}\enskip i\in A_0,~j\in A_1, \enskip\mbox{then}\enskip \left| \Omega \left( \overline{v}_i,%
\overline{\overline{v}}_{i+1}\right) \right| -1+\left| \Omega (\overline{v}%
_j,\overline{\overline{v}}_{j+1})\right| -1\geq 2\mu \left( T\right) . 
$$

Observing that $\left| A_0\right| \geq \left| A_1\right| $, we obtain%
$$
\sum^\kappa_{i=1}(\mid \Omega (\overline{v}_i,\overline{%
\overline{v}}_{i+1})\mid -1)=\sum_{i\in A_0}(\mid \Omega (%
\overline{v}_i,\overline{\overline{v}}_{i+1})\mid -1)+
\sum_{i\in A_1}(\mid \Omega (\overline{v}_i,\overline{\overline{v}}_{i+1})\mid -1)
$$
$$ 
\geq (\left| A_0\right| -\left| A_1\right| )\mu (T)+2\left| A_1\right| \mu
(T)=(\left| A_0\right| +\left| A_1\right| )\mu (T)=\kappa\mu (T). 
$$

{\bf Case b3.2.}$\quad 6\leq h\leq 7.$

Let $i\in A_1.$ Assume w.l.o.g that $\overline{v}_i=u_1,\overline{\overline{v%
}}_{i+1}=u_s$ for some $s\in \{1,...,h\}$. We will write $i\in A_1^{*}$
if and only if $\left| \Omega \left( \overline{v}_i,\overline{\overline{v}}%
_{i+1}\right) \right| -1\geq \beta _j$ for some $j\in \left\{
1,h,s-1,s\right\} $.

{\bf Case b3.2.1.}$\quad A_1=A_1^{*}.$

Let $i\in A_1^{*}$ and let $\overline{v}_i=u_1,\overline{\overline{v}}%
_{i+1}=u_s$~$\left( s\in \{1,...,h\}\right) $. By the definition, $\left|
\Omega \left( \overline{v}_i,\overline{\overline{v}}_{i+1}\right) \right|
-1\geq \beta _j$ for some $j\in \left\{ h,s-1,s\right\} $, say $j=s$. Since $%
6\leq h\leq 7,$ there are at least three pairwise different integers $%
f_1,f_2,f_3$ in $\{1,...,h\}\backslash\left\{ 1,h,s-1\right\} $. By $(i1){\bf ,}%
\left| \Omega \left( \overline{v}_i,\overline{\overline{v}}_{i+1}\right)
\right| -1\geq \max \left( \beta _{f_1},\beta _{f_2},\beta _{f_3}\right) .$
So,%
$$
\mbox{if}\enskip i\in A_1, \enskip\mbox{then}\enskip \left| \Omega \left( \overline{v}_i,\overline{%
\overline{v}}_{i+1}\right) \right| -1\geq \frac 1h(
\sum^h_{i=1}\beta _i-\beta _1-\beta _h-\beta _{s-1}+\beta
_{f_1}+\beta _{f_2}+\beta _{f_3}) 
$$
and hence we can argue exactly as in case $h\geq 8.$

{\bf Case b3.2.2.}$\quad A_1\not =A_1^{*}.$

Let $A_1=\left\{ i,j\right\} $, where $i\not \in A_1^{*}$ and let $\overline{%
v}_i=\overline{v}_j=u_1,~\overline{\overline{v}}_{i+1}=u_s,~\overline{%
\overline{v}}_{j+1}=u_r$ for some $s,r\in \{1,...,h\}$~$\left( s\leq
r\right) $. By $(i8)$ and $(i9){\bf ,~}4\leq s\leq r\leq h-1$. If $s=r$, then
it is easy to see (by definition 3.9) that either $u_1\in \overline{U}_0$ or 
$u_s\in \overline{U}_0$ implying by $(i6)$ that $i\in A_1^{*}$ , a
contradiction. So, assume $s\not =r$, i. e. ${\bf ~}4\leq s<r\leq h-1$.

{\bf Case b3.2.2.1.}$\quad h=7.$

{\bf Case b3.2.2.1.1.}$\quad s=4$ and $r=5.$

By $(i5),$ either $\left| \Omega \left( \overline{v}_i,\overline{\overline{v}%
}_{i+1}\right) \right| -1\geq \beta _4$ or $\left| \Omega \left( \overline{v}%
_j,\overline{\overline{v}}_{j+1}\right) \right| -1\geq \beta _4$. Since $i%
\not \in A_1^{*}$, we have $\left| \Omega \left( \overline{v}_j,\overline{%
\overline{v}}_{j+1}\right) \right| -1\geq \beta _4$. Using $(i1)$, we get
$$
\left| \Omega \left( 
\overline{v}_i,\overline{\overline{v}}_{i+1}\right) \right| -1\geq \frac 1h(%
\sum ^7_{i=1}\beta _i-\beta _1-\beta _7-\beta
_3-\beta _4+2\beta _2+\beta _5+\beta _6),
$$
$$ \left| \Omega \left( \overline{v%
}_j,\overline{\overline{v}}_{j+1}\right) \right| -1\geq \frac 1h(
\sum ^7_{i=1}\beta _i-\beta _1-\beta _7-\beta _5+2\beta _3+\beta
_6). 
$$

Using also all $\kappa-2$ inequalities of type (\ref{17}), we obtain the desired
result as in case $h\geq 8$.

{\bf Case b3.2.2.1.2.}$\quad s=4$ and $r=6.$

By $(i1)$ and $(i9){\bf ,}$%
$$
\begin{array}{l}
\left| \Omega \left( 
\overline{v}_i,\overline{\overline{v}}_{i+1}\right) \right| -1\geq \frac 1h(%
\sum ^7_{i=1}\beta _i-\beta _1-\beta _7-\beta
_3-\beta _4+2\beta _2+\beta _5+\beta _6), \\ \\ \left| \Omega (\overline{v}_j,%
\overline{\overline{v}}_{j+1})\right| -1\geq \frac 1h(
\sum ^7_{i=1}\beta _i-\beta _1-\beta _5+\beta _3+\beta _4). 
\end{array}
$$

Apply the arguments in case b3.2.2.1.1.

{\bf Case b3.2.2.1.3.}$\quad s=5$ and $r=6.$

By $(i1){\bf ,}(i5)$ and $(i9){\bf ,}$%
$$
\begin{array}{l}
\left| \Omega \left( 
\overline{v}_i,\overline{\overline{v}}_{i+1}\right) \right| -1\geq \frac 1h(%
\sum ^7_{i=1}\beta _i-\beta _1-\beta _7-\beta
_4-\beta _5+2\beta _2+\beta _3+\beta _6),\\ \\ \left| \Omega \left( \overline{v%
}_j,\overline{\overline{v}}_{j+1}\right) \right| -1\geq \frac 1h(
\sum ^7_{i=1}\beta _i-\beta _1-\beta _5+\beta _3+\beta _4). 
\end{array}
$$

Apply the arguments in case b3.2.2.1.1.

{\bf Case b3.2.2.2.}$\quad h=6.$

Clearly $s=4,r=5$. By $(i1){\bf ,}(i5)$ and $(i9){\bf ,}$%
$$
\begin{array}{l}
\left| \Omega \left( 
\overline{v}_i,\overline{\overline{v}}_{i+1}\right) \right| -1\geq \frac 1h(%
\sum ^6_{i=1}\beta _i-\beta _1-\beta _6-\beta
_3-\beta _4+2\beta _2+2\beta _5), \\ \left| \Omega (\overline{v}_j,\overline{%
\overline{v}}_{j+1})\right| -1\geq \frac 1h(
\sum ^6_{i=1}\beta _i-\beta _1+\beta _3). 
\end{array}
$$

Apply the arguments in case b3.2.2.1.1.

{\bf Case b3.3.}$\quad h=5.$

Let $A_1=\left\{ i,j\right\} $ and $\overline{v}_i=\overline{v}_j=u_1,~%
\overline{\overline{v}}_{i+1}=u_s,~\overline{\overline{v}}_{j+1}=u_r$ for
some $s,r\in \{1,...,h\}$~$\left( s\leq r\right) $. By $(i8)$ and $(i9)%
{\bf ,}$ $s=r=4$ and we can reach a contradiction as in case b3.2.2.

{\bf Case b4.}$\quad \left| A_{11}\right| =1.$

Let $A_{11}=\left\{ i\right\} $ and $\overline{v}_i=u_1,~\overline{\overline{%
v}}_{i+1}=u_s$ for some $s\in \{1,...,h\}$~.

{\bf Case b4.1.}$\quad h=5.$

By $(i8)$ and $(i9){\bf ,\ }s=4.$ Also, by $(i1)$ and $(i9)$,%
$$
\left| \Omega \left( \overline{v}_i,\overline{\overline{v}}_{i+1}\right)
\right| -1\geq \frac 15(\sum ^5_{i=1}\beta _i-\beta
_1-\beta _3+\beta _2+\beta _4). 
$$

Apply the arguments in case b3.2.2.1.1.

{\bf Case b4.2.}$\quad h\geq 6.$

There are at least two distinct integers $f_1,f_2$ in $\{1,...,h\}%
-\left\{ 1,h,s-1,s\right\} $. By $(i1),$%
$$
\left| \Omega \left( \overline{v}_i,\overline{\overline{v}}_{i+1}\right)
\right| -1\geq \frac 1h(\sum ^h_{i=1}\beta _i-\beta
_1-\beta _h-\beta _{s-1}-\beta _s+2\beta _{f_1}+2\beta _{f_2}). 
$$

Since $\mid A_0\mid \geq 4-\mid A_1\mid =3$, we have at least two
inequalities of type (\ref{17}). So, we can argue as in case b3.2.2.1.1.\quad$%
\Delta $

By claims 5 and 6, $c\geq \sum_{i=1}^r\varphi _{v_i^{*}}+2\kappa+\kappa\mu \left(
T\right) $ and the result follows as in case 1.1.1.1.

{\bf Case 2.}$\quad \kappa\geq 4,h\leq 4.$

Since $h\geq \kappa$, we have $h=\kappa=4$. Then
there are four vertex-disjoint $(H,C)-$paths. It can be easily cheeked that $%
c\geq 18$. If $\delta \leq 6$, then $c\geq 18\geq 20\left( \delta +2\right)
/9=\left( h+1\right) \kappa(\delta +2)/(h+\kappa+1)$. Let $\delta \geq 7$. Using $(d3)$
we can show that $\varphi _i\leq 3$ for some $i\in \{1,2,3,4\}$, i.e. $%
\max_i\psi _i\geq \delta -3$. Then by lemma 3, $c\geq
\sum_{i=1}^4\left( \delta -\varphi _i\right) +\delta -3=5\delta
-3-\sum_{i=1}^4\varphi _i.$ If $\sum_{i=1}^4\varphi _i\leq 12$, then $c\geq
5\delta -15\geq \left( h+1\right) \kappa\left( \delta +2\right) /\left(
h+\kappa+1\right) $. So, it sufficies to prove $\sum_{i=1}^4\varphi _i\leq 12$.

{\bf Case 2.1.}$\quad $Either $\left| U_0\right| =0$ or $\left| U_0\right|
=4.$

It follows that $\varphi _i\leq 3$ for each $i\in \{1,2,3,4\}$.

{\bf Case 2.2}.$\quad \left| U_0\right| =3.$

Assume w.l.o.g. that $\overline{U}_0=\left\{ u_1\right\} $. If $u_3\stackrel{%
o}{u}_1\not \in E$, then it is easy to see that $\varphi _i\leq 3$ for each $%
i\in \{1,2,3,4\}$. Otherwise $(u_3\stackrel{o}{u}_1\in E),~u_2u_4\not \in
E$ and hence $\varphi _1\leq 3,\varphi _3\leq 4,\varphi _2\leq 2$ and $%
\varphi _4\leq 2$.

{\bf Case 2.3.}$\quad \left| U_0\right| =2.$

By symmetry, we can distinguish the following two cases.

{\bf Case 2.3.1.}$\quad \overline{U}_0=\left\{ u_1,u_4\right\} .$

If $u_3\stackrel{0}{u}_1\not \in E$ and $u_2\stackrel{0}{u}_4\not \in E$,
then clearly $\varphi _i\leq 3$ for each $i\in \{1,2,3,4\}$. Assume
w.l.o.g. that $u_3\stackrel{0}{u}_1\in E.$ We can assume also $u_2\stackrel{0%
}{u}_4\not \in E$, since otherwise the cycle $u_1\stackrel{0}{u}_1u_3u_4%
\stackrel{0}{u}_4u_2u_1$ is larger than $H$, which is impossible. Then
clearly $\varphi _1\leq 3,\varphi _4\leq 3,\varphi _3\leq 4$ and $\varphi
_2\leq 2$.

{\bf Case 2.3.2.}$\quad \overline{U}_0=\left\{ u_1,u_3\right\} .$

It is easy to see that $\varphi _i\leq 3$ for each $i\in \{1,2,3,4\}$.

{\bf Case 2.4.}$\quad \left| U_0\right| =1.$ 

Returning to the proof of lemma
3, we can see that in this special case the lower bound in lemma 3 can be
improved by a unit. So, it suffices to show $\sum_{i=1}^4\varphi _i\leq 13$.
Denoting $U_0=\left\{ u_4\right\} $, we see that $\varphi _1\leq 3,\varphi
_2\leq 3,\varphi _3\leq 3,\varphi _4\leq 4$ and the result holds immediately.

{\bf Case 3.}$\quad \kappa\leq 3.$\\

{\bf Claim 7.}\quad Let $\kappa\in \{2,3\}$ and $h\geq \kappa$. If there are no $\kappa+1$ vertex-disjoint 
$(H,C)-$paths, then $c\geq \min \left( \kappa\left( h+1\right) ,\kappa\left( \delta -\kappa+4\right) \right).$

{\bf Proof of claim 7.}\quad {\bf Case d1.}\quad $\kappa=3.$

Assume w.l.o.g. that $E_1,E_2$ and $E_3$ are $T-$transformed. We now prove
that $\mid w_1\overrightarrow{C}w_2\mid -1\geq \min \left( h+1,\delta
+1\right) $. If $v_2=v_1^{+}$, then clearly 
$$
\mid w_1\overrightarrow{C}w_2\mid -1\geq \mid w_1E_1v_1\overleftarrow{H}%
v_2E_2w_2\mid -1\geq h+1. 
$$
Let $v_2\not =v_1^{+}.$ Walking along $\overrightarrow{H}$ from $v_1$ to $%
v_2^{-}$ we stop at the first vertex $z$ with either $\stackrel{\wedge }{z}%
w_2\in E$ or $\stackrel{\wedge }{z}w_1\not \in E$ or $z=v_2^{-}$. If $%
\stackrel{\wedge }{z}w_2\in E$ or $z=$ $v_2^{-}$, then clearly $\mid w_1%
\overrightarrow{C}w_2\mid -1\geq h+1$. Let $\stackrel{\wedge }{z}w_2\in E$
and $\stackrel{\wedge }{z}w_1\not \in E$. If $\stackrel{\wedge }{z}w\in E$
for some $w\in V\left( C\right) \backslash\left\{ w_1,w_2,w_3\right\} $, then there
are $4$ vertex-disjoint $(H,C)-$paths, contradicting our assumption. So, $N(%
\stackrel{\wedge }{z})\bigcap V\left( C\right) \subseteq \left\{ w_3\right\} 
$, i.e. $\varphi _z\geq \delta -1$ and $h\geq \varphi _z+1\geq \delta $. By $%
(g6)$, $\left| O\left( z^{-},v_2\right) \right| -1\geq \gamma _z\geq $ $%
\varphi _z\geq \delta -1$ implying that $\mid w_1\overrightarrow{C}w_2\mid
-1\geq \delta +2$. Thus we have proved $\mid w_1\overrightarrow{C}w_2\mid
-1\geq \min \left( h+1,\delta +1\right) .$ By symmetry, we have similar
inequalities for segments $w_2\overrightarrow{C}w_3$ and $w_3\overrightarrow{%
C}w_1$ and the result holds from $h+1\geq \delta +1$.

{\bf Case d2.}\quad $\kappa=2.$

Apply the arguments in case 1. Claim 7 is proved.\quad$\Delta $

{\bf Case 3.1.}\quad $\kappa=3.$

We can assume that there are no 4 vertex-disjoint $(H,C)-$paths, since
otherwise%
$$
c\geq \frac{\left( h+1\right) 4}{h+4+1}\left( \delta +2\right) >\frac{\left(
h+1\right) \kappa}{h+\kappa+1}\left( \delta +2\right) . 
$$
Then by claim 7 we can distinguish the following two cases.

{\bf Case 3.1.1.}\quad $c\geq 3\left( h+1\right) .$

If $h\geq \delta -2$, then $c\geq 3\left( h+1\right) \geq 3\left( h+1\right)
\left( \delta +2\right) /\left( h+4\right) $. Otherwise, the result holds
from $c\geq 3\left( \delta -1\right) $ (see \cite{12}).

{\bf Case 3.1.2.}\quad $c\geq 3\left( \delta +1\right) .$

If $h\leq 3\delta +2$, then $c\geq 3\left( \delta +1\right) \geq 3\left(
h+1\right) \left( \delta +2\right) /\left( h+4\right) $. Let $h\geq 3\left(
\delta +1\right) $. Observing that $c\geq 3\left( h/2+2\right) $ (by
standard arguments) we obtain the result immediately.

{\bf Case 3.2.}\quad $\kappa=2.$

Apply the arguments in case 1.

Thus we have proved the theorem for $h$ the length of a longest cycle in $G\backslash C$. Observing
that%
$$
c\geq \frac{\left( h+1\right) \kappa}{h+\kappa+1}\left( \delta +2\right) >\frac{\left(
h^{^{\prime }}+1\right) \kappa}{h^{^{\prime }}+\kappa+1}\left( \delta +2\right) 
$$
for any $h^{^{\prime }}<h$ , we complete the proof of the theorem.\quad$\Delta $

\section*{Acknowledgment}

The author is very grateful to N.\ K. Khachatrian for his careful reading
and his corrections for the manuscript.


\begin{thebibliography}{99}
\bibitem{1}  J. A. Bondy and U. S. R. Murty, Graph Theory with Applications.
Macmillan, London and Elsevier, New York (1976).

\bibitem{2}  G. A. Dirac, Some theorems on abstract graphs, Proc. London,
Math. Soc., 2 (1952), 69--81.

\bibitem{3}  Zh.G.Nikoghosyan, On maximal cycle of a graph, DAN Arm.SSR v.LXXII
2 (1981) 82-87 (in Russian).

\bibitem{4}  Zh. G. Nikoghosyan, Path-Extensions and Long Cycles in Graphs,
Transactions of the Institute for Informatics and Automation Problems of the
NAS of RA and of the Yerevan State University, Mathematical Problems of
Computer Science 19 (1998), 25-31.

\bibitem{5}  Zh. G. Nikoghosyan, Cycle-Extensions and Long Cycles in Graphs,
Transactions of the Institute for Informatics and Automation Problems of the
NAS of RA and of the Yerevan State University, Mathematical Problems of
Computer Science 21 (2000), 121-128.

\bibitem{6}	H.-J. Voss and C. Zuluaga, Maximale gerade und ungerade Kreise in Graphen I, Wiss. Z. Tech. Hochschule Ilmenau, 23 (1977) 57-70.

\end{thebibliography}
\end{document}